\documentclass[11pt]{amsart}
\usepackage{amsmath, amssymb, amsthm, amsfonts}
\usepackage{hyperref,doi}
\usepackage{epsf}
\usepackage{enumerate}
\usepackage{graphicx}
\usepackage{array}
\usepackage[margin=1in]{geometry}
\newcommand{\ZZ}{\ensuremath{\mathbb{Z}}}

\newcommand{\C}{\mathbb{C}}
\newcommand{\R}{\mathbb{R}}

\newcommand{\W}{\mathfrak{W}}

\newcommand{\sym}{\mathfrak{S}}

\newcommand{\G}{\mathcal{G}}

\newcommand{\el}{\mathcal{L}}

\newcommand{\p}{\mathcal{P}}

\newcommand{\I}{\mathcal{I}}

\newcommand{\M}{\mathcal{M}}

\newcommand{\khsr}{\Delta_{n,k}^{stab(r)}}
\newcommand{\rstablenk}{\Delta_{n,k}^{stab(r)}}
\newcommand{\rstablentwo}{\Delta_{n,2}^{stab(r)}}
\newcommand{\rstablentwoint}{\Delta_{n,2}^{stab(r)^\circ}}
\newcommand{\hypernk}{\Delta_{n,k}}
\newcommand{\hyperntwo}{\Delta_{n,2}}
\newcommand{\twostablentwo}{\Delta_{n,2}^{stab(2)}}
\newcommand{\threestablentwo}{\Delta_{n,2}^{stab(3)}}

\newcommand{\cd}[2]{\operatorname{cd}({#1},{#2})}
\newcommand{\arc}[2]{\operatorname{arc}({#1},{#2})}
\newcommand{\sort}[1]{\operatorname{sort}({#1})}
\newcommand{\U}[2]{\operatorname{\emph{U}}({#1},{#2})}
\newcommand{\V}[2]{\operatorname{\emph{V}}({#1},{#2})}
\newcommand{\adjvar}[2]{\operatorname{adj_{{#1}}({#2})}}
\newcommand{\adj}[1]{\operatorname{adj_{r}({#1})}}
\newcommand{\Adj}[1]{\operatorname{Adj_{\emph{r}}[{#1}]}}
\newcommand{\Adjvar}[2]{\operatorname{Adj_{{#1}}[{#2}]}}
\newcommand{\hstar}[1]{\operatorname{h^\ast\left({#1};x\right)}}
\newcommand{\hstarint}[1]{\operatorname{h^\ast\left({#1};\frac{1}{x}\right)}}

\newcommand{\floorntwo}{\left\lfloor\frac{n}{2}\right\rfloor}
\newcommand{\ceilingntwo}{\left\lceil\frac{n}{2}\right\rceil}

\newcommand{\trink}{\nabla_{n,k}}
\newcommand{\trinkr}{\nabla_{n,k}^r}
\newcommand{\trintwo}{\nabla_{n,2}}
\newcommand{\trintwor}{\nabla_{n,2}^r}
\newcommand{\trintworplus}{\nabla_{n,2}^{r+1}}

\newtheorem{thm}{Theorem}[section]

\newtheorem{quest}[thm]{Question}
\newtheorem{conj}[thm]{Conjecture}
\theoremstyle{definition}

\newtheorem{defn}[thm]{Definition}
\newtheorem{cor}[thm]{Corollary}
\newtheorem{lem}[thm]{Lemma}
\newtheorem{rmk}[thm]{Remark}
\newtheorem{example}[thm]{Example}
\newtheorem{prop}[thm]{Proposition}

\DeclareMathOperator{\vol}{vol}

\DeclareMathOperator{\Des}{Des}
\DeclareMathOperator{\colex}{colex}
\DeclareMathOperator{\adjr}{adj_r}

\begin{document}

\begin{abstract}
Hypersimplices are well-studied objects in combinatorics, optimization, and representation theory.
For each hypersimplex, we define a new family of subpolytopes, called $r$-stable hypersimplices, and show that a well-known regular unimodular triangulation of the hypersimplex restricts to a triangulation of each $r$-stable hypersimplex.
For the case of the second hypersimplex defined by the two-element subsets of an $n$-set, we provide a shelling of this triangulation that sequentially shells each $r$-stable sub-hypersimplex.
In this case, we utilize the shelling to compute the Ehrhart $h^\ast$-polynomials of these polytopes, and the hypersimplex, via independence polynomials of graphs.  
For one such $r$-stable hypersimplex, this computation yields a connection to CR mappings of Lens spaces via Ehrhart-MacDonald reciprocity.
\end{abstract}

\title{Shellability, Ehrhart Theory, and $r$-stable Hypersimplices}

\author{Benjamin Braun}
\address{Department of Mathematics\\
         University of Kentucky\\
         Lexington, KY 40506--0027}
\email{benjamin.braun@uky.edu}

\author{Liam Solus}
\address{Institute of Science and Technology Austria \\ 
	Klosterneuburg, Austria \\ 
	and \\
	Institute for Data, Systems, and Society\\
         Massachusetts Institute of Technology \\
         Cambridge, MA, USA}
\email{lsolus@mit.edu}

\date{16 March 2016}

\thanks{Benjamin Braun was partially supported by the National Security Agency through award H98230-13-1-0240.  
Liam Solus was partially supported by a 2014 National Science Foundation/Japan Society for the Promotion of Science East Asia and Pacific Summer Institute Fellowship.}


%52B05 = Polytopes and Polyhedra - Combinatorial Properties (number of faces, shortest paths, etc...)
%52B20 = Polytopes and Polyhedra - Lattice Polytopes (including connections with commutative algebra and algebraic geometry)


\maketitle

\section{Introduction}
Fix integers $0<k<n$.  
We let $[n]:=\{1,2,\ldots,n\}$ and let $[n]\choose k$ denote the collection of $k$-element subsets  of $[n]$.  
The \emph{characteristic vector} of a $k$-subset $I$ of $[n]$ is the $(0,1)$-vector $\epsilon_I:=(\epsilon_1,\ldots,\epsilon_n)$ such that $\epsilon_i=1$ for $i\in I$ and $\epsilon_i=0$ for $i\notin I$.  
The \emph{$(n,k)$-hypersimplex}, denoted $\Delta_{n,k}$, is the $(n-1)$-dimensional polytope that is the convex hull of the characteristic vectors of all $k$-subsets of $[n]$.  
That is, $\Delta_{n,k}$ is the convex hull of all $(0,1)$-vectors in $\R^n$ containing precisely $k$ nonzero terms.  
Hypersimplices appear naturally in algebraic and geometric contexts, as well as in pure and applied combinatorial contexts.  
In \cite{stan2}, Stanley geometrically proved that the volume of the $(n,k)$-hypersimplex is the Eulerian number $A_{k,n-1}$.  
De Loera, Sturmfels, and Thomas studied the connection between triangulations of the hypersimplex and Gr\"{o}bner bases via toric algebra \cite{del2}.  
In \cite{lam}, Lam and Postnikov showed four useful triangulations of the hypersimplex are all identical.  
In \cite{kat}, Katzman gave an algebraic description of the Ehrhart $h^\ast$-vector of the hypersimplex, and in \cite{li}, Li gave a second interpretation in terms of exceedences and descents.  
These many investigations have proven fruitful for our understanding of these fundamental polytopes, but many interesting questions about $\Delta_{n,k}$ still remain unanswered.

Many of the unanswered questions pertaining to the hypersimplex lie in the field of Ehrhart theory.  
A \emph{lattice polytope} $\p$ of dimension $d$ is the convex hull in $\R^n$ of finitely many points in $\ZZ^n$ that together affinely span a $d$-dimensional hyperplane in $\R^n$.  For $t\in\ZZ_{>0}$, set $t\p:=\{tp:p\in\p\}$, and let $L_\p(t)=|\ZZ^n\cap t\p|$.  In \cite{ehr}, Ehrhart proved that with the polynomial basis $\{{t+d-i\choose d}:i\in[0,d]\cap\ZZ\}$, for any lattice polytope $\p$ we have
\begin{equation*}
L_\p(t)=\sum_{i=0}^dh_i^*{t+d-i\choose d}.
\end{equation*}
Stanley then proved that the coefficients $h_i^*$ are all nonnegative integers \cite{stan3}.  The polynomial $L_\p(t)$ is called the \emph{Ehrhart polynomial of $\p$} and has connections to commutative algebra, algebraic geometry, combinatorics, and discrete and convex geometry.  However, these polynomials are not-well understood in many cases.  Given the above polynomial representation of $L_\p(t)$ it is common to study the coefficient vector $h^*(\p):=(h_1^*,\ldots,h_d^*)$, which is often referred to as the \emph{h-star vector} or \emph{$\delta$-vector} of $\p$.  This vector encodes useful information about the polytope $\p$.  For example, $\vol(\p)=\dfrac{\sum_ih_i^*}{d!}$, where $\vol(\p)$ denotes the Euclidean Volume (Lebesgue measure) of $\p$ with respect to the integer lattice contained in the hyperplane spanned by $\p$.  We call the sum $\sum_ih_i^*$ the \emph{normalized volume of $\p$}.  Another useful (but not always achieved) property of these vectors is \emph{unimodality}.  A vector $(x_0,x_1,\ldots,x_d)$ is called \emph{unimodal} if there exists an index $j$, $0\leq j\leq d$, such that $x_{i-1}\leq x_i$ for $i\leq j$, and $x_i\geq x_{i+1}$ for $i\geq j$.  Unimodality of the $h^*$-vector has interesting algebraic implications and consequently is a widely sought property of these vectors \cite{hibi,kat,stan}.
In the special case of the hypersimplex, Haws, De Loera, and K\"{o}ppe computationally verified that the $h^\ast$-vector of $\hypernk$ is unimodal for $n$ as large as $75$ \cite{del}.  Since the $h^\ast$-vector of $\hypernk$ is not symmetric, these results are intriguing from both the perspective of investigation of  the $h^\ast$-vectors of hypersimplices as well as the perspective of investigation of unimodal sequences.

The main purpose of this article is to introduce a new variation on the  hypersimplices and examine how this variation can help us to better understand their $h^\ast$-vectors.
In the following, for each integer $r$ satisfying $1\leq r\leq\left\lfloor\frac{n}{k}\right\rfloor$ we identify a subpolytope of the $(n,k)$-hypersimplex that we call the \emph{$r$-stable $(n,k)$-hypersimplex}.
These subpolytopes are constructed in a manner that forms a nested chain within the $(n,k)$-hypersimplex, i.e. the $r$-stable $(n,k)$-hypersimplex is a subpolytope of the $(r-1)$-stable $(n,k)$-hypersimplex. 
This nesting is the seed of a useful geometric relationship between these polytopes.
In section \ref{triangulation} we describe the nature of this geometric relationship; our main result from this section is the following.

\begin{thm}\label{inducedtriangulationthm}
The {\bf circuit triangulation}, a well-known regular unimodular triangulation of the $(n,k)$-hypersimplex, restricts to a triangulation of each $r$-stable $(n,k)$-hypersimplex.
\end{thm}

The circuit triangulation of $\hypernk$ is first defined in \cite{lam}, and will serve as one of our most important computational tools.  
The power of Theorem \ref{inducedtriangulationthm} lies in the unimodularity of the induced triangulations.  
A result of Stanley \cite{stan3} says that a shellable unimodular triangulation of an integral polytope may be used to compute the $h^\ast$-polynomial of the polytope.  
Consequently, if there exists a shelling of the circuit triangulation that proceeds with respect to the nesting of $r$-stable hypersimplices then we may inductively compare these $h^\ast$-polynomials.
In section \ref{shelling} we demonstrate that such a shelling exists in the case when $k=2$.

\begin{thm}\label{introshellingthm}
There exists a shelling of the circuit triangulation that first builds the $r$-stable $(n,2)$-hypersimplex and then builds the $(r-1)$-stable $(n,2)$-hypersimplex for every $1\leq r <\floorntwo$.  
\end{thm}
\noindent In section~\ref{extending results}, we introduce a conjectured method for extending this shelling to $r$-stable $(n,k)$-hypersimplices for arbitrary $k>2$.  

In section \ref{hstarpolynomials}, we utilize the shelling of Theorem \ref{introshellingthm} to consecutively compute the $h^\ast$-polynomials of the $r$-stable $(n,2)$-hypersimplices, and the $(n,2)$-hypersimplex containing them.  
The main result of this section is that the $h^\ast$-polynomial of all these polytopes, including the hypersimplex, may be computed via sums of independence polynomials of graphs.

\begin{thm}\label{introsumsofindependencepolysthm} 
For each $1\leq j\leq\floorntwo$ and $\ell\in[n]$, there exist graphs $G_{n,j,\ell}$ such that the $h^\ast$-polynomial of the $r$-stable $(n,2)$-hypersimplex equals
	$$1+x\left(\sum_{j=r}^{\floorntwo-1}\sum_{\ell=1}^nI\left(G_{n,j,\ell};x\right)\right),$$
where $I\left(G_{n,j,\ell};x\right)$ denotes the independence polynomial of $G_{n,j,\ell}$.
\end{thm}

From these computations we discover that the $h^\ast$-polynomials of the $r$-stable hypersimplices are fascinating in their own right.  
In particular, with Theorem \ref{lucasthm} we will see that for each odd $n$, the $n^{th}$ Lucas polynomials arise as the $h^\ast$-polynomials of a collection of $r$-stable hypersimplices.  
A consequence of Theorem \ref{lucasthm} is a connection, via Ehrhart-MacDonald reciprocity, between these $r$-stable hypersimplices and CR mappings of Lens spaces into complex unit spheres.  

We end section \ref{hstarpolynomials} with a discussion of unimodality.  
It is known that the $h^\ast$-polynomial of the second hypersimplex is unimodal \cite{kat}, and we demonstrate that this is also true for a collection of $r$-stable hypersimplices within.  
Collectively, these results suggest that the geometric relationship between the $r$-stable hypersimplices and the hypersimplex containing them can be useful from an Ehrhart--theoretical perspective, while also suggesting that the $r$-stable hypersimplices have interesting structure in their own right.

%%%%%%%%%%%%%%%%%%%%%%%%%%%%%%%%%%%%%%%%%%%%%%%%%%%%%%%%%%%%%%%
%%%%%%%%%%%%%%%%%%%%%%%%%%%%%%%%%%%%%%%%%%%%%%%%%%%%%%%%%%%%%%%
%%%%%%%%%%%%%%%%%%%%%%%%%%%%%%%%%%%%%%%%%%%%%%%%%%%%%%%%%%%%%%%

\section{The $r$-stable $(n,k)$-hypersimplex}\label{triangulation}
Label the vertices of a regular $n$-gon embedded in $\R^2$ in a clockwise fashion from $1$ to $n$.  We define the \emph{circular distance} between two elements $i$ and $j$ of $[n]$, denoted $\cd{i}{j}$, to be the number of edges in the shortest path between the vertices $i$ and $j$ of the $n$-gon.  We also denote the path of shortest length from $i$ to $j$ by $\arc{i}{j}$.  A subset $S\subset[n]$ is called \emph{$r$-stable} if each pair $i,j\in S$ satisfies $\cd{i}{j}\geq r$.  The \emph{$r$-stable $(n,k)$-hypersimplex}, denoted $\Delta_{n,k}^{stab(r)}$, is the convex hull of the characteristic vectors of all $r$-stable $k$-subsets of $[n]$.  For fixed $n$ and $k$, these polytopes form the nested chain
\begin{equation*}
\Delta_{n,k}\supset\Delta_{n,k}^{stab(2)}\supset\Delta_{n,k}^{stab(3)}\supset\cdots\supset\Delta_{n,k}^{stab\left(\left\lfloor\frac{n}{k}\right\rfloor\right)}.
\end{equation*}

\subsection{A well-studied triangulation of the hypersimplex}

In \cite{lam}, Lam and Postnikov compare four different triangulations of the hypersimplex, and show that they are identical.  
While these triangulations possess the same geometric structure the constructions are all quite different, and consequently each one reveals different information about the triangulation's geometry.
Here, we utilize properties of two of these four constructions.
The first is a construction given by Sturmfels in \cite{stur} using techniques from toric algebra.  
The second construction, known as the \emph{circuit triangulation}, is introduced in \cite{lam} by Lam and Postnikov.
We will show that this triangulation restricts to a triangulation of the $r$-stable hypersimplex.

\subsubsection{Sturmfels' Triangulation}
We recall the description of this triangulation presented in \cite{lam}.
Let $I$ and $J$ be two $k$-subsets of $[n]$ and consider their multi-union $I\cup J$.  Let $\sort{I\cup J}=(a_1,a_2,\ldots,a_{2k})$ be the unique nondecreasing sequence obtained by ordering the elements of the multiset $I\cup J$ from least-to-greatest.  
Now let $\U{I}{J}:=\{a_1,a_3,\ldots,a_{2k-1}\}$ and $\V{I}{J}:=\{a_2,a_4,\ldots,a_{2k}\}$.  
As an example consider the $4$-subsets of $[8]$, $I=\{1,3,4,6\}$ and $J=\{3,5,7,8\}$.  
For this pair of subsets we have that $\sort{I\cup J}=(1,3,3,4,5,6,7,8)$, $\U{I}{J}=\{1,3,5,7\}$, and $\V{I}{J}=\{3,4,6,8\}$.
The ordered pair of $k$-subsets $(I,J)$ is said to be \emph{sorted} if $I=\U{I}{J}$ and $J=\V{I}{J}$.  
Moreover, an ordered $d$-collection $\I=(I_1,I_2,\ldots,I_d)$ of $k$-subsets  is called \emph{sorted} if each pair $(I_i,I_j)$ is sorted for all $1\leq i<j\leq d$.  
For a sorted $d$-collection $\I$ we let $\sigma_\I$ denote the $(d-1)$-dimensional simplex with vertices $\epsilon_{I_1},\epsilon_{I_2},\ldots,\epsilon_{I_d}$.

\begin{thm}\cite[Sturmfels]{stur}\label{stur1}
The collection of simplices $\sigma_\I$, where $\I$ varies over the sorted collections of $k$-element subsets of $[n]$, forms a triangulation of $\Delta_{n,k}$.  
\end{thm}

Notice that the maximal simplices in this triangulation correspond to the maximal-by-inclusion sorted collections, which all have $d=n$.

This triangulation of $\Delta_{n,k}$ was identified by Sturmfels' via the correspondence between Gr\"{o}bner bases for the toric ideal associated to $\Delta_{n,k}$ and regular triangulations of $\Delta_{n,k}$.  
To construct the toric ideal for $\Delta_{n,k}$ let $k[x_I]$ denote the polynomial ring in the $n \choose k$ variables $x_I$ labeled by the $k$-subsets of $[n]$, and define the semigroup algebra homomorphism
$$\varphi:k[x_I]\longrightarrow k[z_1,z_2,\ldots,z_n]; \qquad \varphi: x_I\longmapsto z_{i_1}z_{i_2}\cdots z_{i_k}, \quad \mbox{for $I=\{i_1,i_2,\ldots, i_k\}$.}$$
The kernel of this homomorphism, $\ker\varphi$, is the toric ideal of $\Delta_{n,k}$.
The correspondence between Gr\"{o}bner bases for $\ker\varphi$ and regular triangulations of $\Delta_{n,k}$ is given as follows.
Any sufficiently generic height vector induces a regular triangulation of $\Delta_{n,k}$.  
On the other hand, such a height vector induces a term order $<$ on the monomials in the polynomial ring $k[x_I]$.  Thus, we may identify a Gr\"{o}bner basis for $\ker\varphi$ with respect to this term order, say $G_<$. 
Moreover, the initial ideal associated to a Gr\"{o}bner basis is square-free if and only if the corresponding regular triangulation is unimodular.  
The details of this correspondence are outlined nicely in \cite{stur}.

\begin{thm}\cite[Sturmfels]{stur}\label{stur2}
The set of quadratic binomials
$$\G_<:=\left\{\underline{x_Ix_J}-x_{\U{I}{J}}x_{\V{I}{J}} : I,J\in{[n]\choose k}\right\}$$
is a Gr\"{o}bner basis for $\ker\varphi$ under some term order $<$ on $k[x_I]$ such that the underlined term is the initial monomial.  In particular, the initial ideal of $G_<$ is square-free, and the simplices of the corresponding unimodular triangulation are $\sigma_\I$, where $\I$ varies over the sorted collections of $k$-element subsets of $[n]$.
\end{thm}

We denote this triangulation of $\Delta_{n,k}$ by $\nabla_{n,k}$, and we let $\max\nabla_{n,k}$ denote the set of maximal simplices in $\nabla_{n,k}$.  In \cite{lam}, Lam and Postnikov prove a more general version of Theorem \ref{stur2} which we will utilize to show that this triangulation restricts to a triangulation of the $r$-stable hypersimplex $\rstablenk$.

\subsubsection{The Circuit Triangulation}   
The second construction of this triangulation that we will utilize first appeared in \cite{lam}, and it arises from examining minimal length circuits in a particular directed graph with labeled edges.
We construct this directed graph as follows.  
Let $G_{n,k}$ be the directed graph with vertices $\epsilon_I$, where $I$ varies over all $k$-subsets of $[n]$.  
For a vertex $\epsilon=(\epsilon_1,\ldots,\epsilon_n)$ of $G_{n,k}$ we think of the coordinate indices $i$ as elements of the cyclic group $\ZZ/n\ZZ$.  
Hence, $\epsilon_{n+1}=\epsilon_1$.  
We construct the directed, labeled edges of $G_{n,k}$ as follows.
Suppose $\epsilon=(\epsilon_1, \ldots, \epsilon_n)$ and $\epsilon^\prime$ are vertices of $G_{n,k}$ for which $(\epsilon_i,\epsilon_{i+1})=(1,0)$ and the vector $\epsilon^\prime$ is obtained from $\epsilon$ by switching $\epsilon_i$ and $\epsilon_{i+1}$.  
Then we include the directed labeled edge $\epsilon\overset{i}\rightarrow \epsilon^\prime$ in $G_{n,k}$.
Hence, each edge of $G_{n,k}$ is given by shifting a 1 in a vertex $\epsilon$ exactly one entry to the right (modulo $n$), and this can happen if and only if the next place is occupied by a 0.

\begin{figure}
	\centering
	\includegraphics{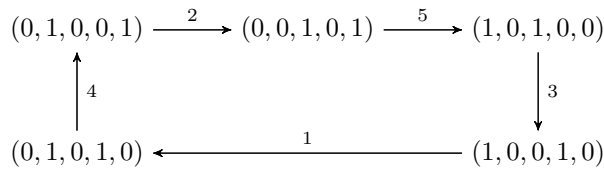}
	\caption{Here is a minimal circuit in the directed graph $G_{5,2}$.}
	\label{fig:circuitexample}
\end{figure}

We are interested in the circuits of minimal possible length in the graph $G_{n,k}$.  
We will call such a circuit \emph{minimal}.  
A minimal circuit in $G_{n,k}$ containing the vertex $\epsilon$ is given by a sequence of edges moving each 1 in $\epsilon$ into the position of the 1 directly to its right.  
Hence, the length of such a circuit is precisely $n$.   
An example of a minimal circuit is given in Figure \ref{fig:circuitexample}.

For a fixed initial vertex, the sequence of labels of edges in a minimal circuit forms a permutation $\omega=\omega_1\omega_2\cdots\omega_n\in \sym_n$, the symmetric group on $n$ elements.  
There is one such permutation for each choice of initial vertex in the minimal circuit.
Hence, a minimal circuit in $G_{n,k}$ corresponds to an equivalence class of permutations in $\sym_n$ where permutations are equivalent modulo cyclic shifts $\omega_1\cdots\omega_n\sim\omega_n\omega_1\cdots\omega_{n-1}$.
In the following, we choose the representative $\omega$ of the class of permutations associated to the minimal circuit for which $\omega_n=n$.
We remark that this corresponds to picking the initial vertex of the minimal circuit to be the lexicographically maximal $(0,1)$-vector in the circuit.  
For example, the lexicographic ordering on the $(0,1)$-vectors in the circuit depicted in Figure \ref{fig:circuitexample} is
$$(1,0,1,0,0)>(1,0,0,1,0)>(0,1,0,1,0)>(0,1,0,0,1)>(0,0,1,0,1),$$

\noindent and the permutation given by reading the edge labels of this circuit beginning at vertex $(1,0,1,0,0)$ is $\omega=31425$.  Thus, we see that $\omega_n=n$ as desired.

\begin{thm}\cite[Lam and Postnikov]{lam}\label{lam1}
A minimal circuit in the graph $G_{n,k}$ corresponds uniquely to a permutation $\omega\in \sym_n$ modulo cyclic shifts.  Moreover, a permutation $\omega\in \sym_n$ with $\omega_n=n$ corresponds to a minimal circuit in $G_{n,k}$ if and only if the inverse permutation $\omega^{-1}$ has exactly $k-1$ descents.  
\end{thm}

We label the minimal circuit in the graph $G_{n,k}$ corresponding to the permutation $\omega\in \sym_n$ with $\omega_n=n$ by $(\omega)$.  Let $v_{(\omega)}$ denote the set of all vertices $\epsilon_I$ of $\Delta_{n,k}$ used by the circuit $(\omega)$, and let $\sigma_{(\omega)}$ denote the convex hull of $v_{(\omega)}$.  

\begin{thm}\cite[Lam and Postnikov]{lam}\label{lam2}
The collection of simplices $\sigma_{(\omega)}$ corresponding to all minimal circuits in $G_{n,k}$ forms the collection of maximal simplices of a triangulation of the hypersimplex $\Delta_{n,k}$.  This triangulation is identical to the triangulation $\nabla_{n,k}$.  
\end{thm}

\noindent We call this construction of $\nabla_{n,k}$ the \emph{circuit triangulation}.  
To simplify notation we will often write $\omega$ for the simplex $\sigma_{(\omega)}\in\nabla_{n,k}$.

\subsection{The induced triangulation of the $r$-stable $(n,k)$-hypersimplex}
Let $\M$ be a collection of $k$-subsets of $[n]$, and let $\p_\M$ denote the convex hull in $\R^n$ of the $(0,1)$-vectors $\{\epsilon_I:I\in\M\}$.  Notice that $\p_\M$ is a subpolytope of $\hypernk$.  The collection $\M$ is said to be \emph{sort-closed} if for every pair of subsets $I$ and $J$ in $\M$ the subsets $\U{I}{J}$ and $\V{I}{J}$ are also in $\M$.  
In \cite{lam}, Lam and Postnikov proved the following theorem.

\begin{thm}\cite[Lam and Postnikov]{lam}\label{lam4}
The triangulation $\nabla_{n,k}$ of the hypersimplex $\Delta_{n,k}$ induces a triangulation of the polytope $\p_\M$ if and only if $\M$ is sort-closed.  
\end{thm}

Using this theorem, we can provide a proof of Theorem~\ref{inducedtriangulationthm}.

\subsubsection{Proof of Theorem~\ref{inducedtriangulationthm}.}
\label{stablesortclosed}
Fix an integer $0<r\leq\left\lfloor\frac{n}{k}\right\rfloor$.  
Let $\M$ be the collection of $r$-stable $k$-subsets of $[n]$.  
We now show that the triangulation $\nabla_{n,k}$ induces a triangulation of the $r$-stable $(n,k)$-hypersimplex $\p_\M=\rstablenk$.
By Theorem \ref{lam4}, it suffices to show that the collection $\M$ is sort-closed.  Let $I$ and $J$ be two elements of $\M$, and consider $\sort{I\cup J}=(a_1,a_2,\ldots,a_{2k})$.  Suppose for the sake of contradiction that for some $i$, $a_{i+2}=a_i+t$ for some $t\in[r-1]$.  Here we think of our indices and addition modulo $n$.  We remark that $t$ must be nonzero since the multiplicity of each element of $[n]$ in $\sort{I\cup J}$ is at most two.  Without loss of generality, we assume that $a_i\in I$.  Hence, $a_{i+2}\in J$ since $I$ is $r$-stable and $\cd{a_i}{a_{i+2}}<r$.  Since $\sort{I\cup J}$ is nondecreasing it follows that $a_{i+1}=a_i+j$ for some $j\in\{0,1,\ldots,t\}$.  First consider the cases where $j=0$ and $j=t$.  In the former case we have that $a_{i+1}=a_i$, and in the latter case $a_{i+1}=a_{i+2}$.  Hence, in the former case, the multiplicity of $a_{i}$ in $\sort{I\cup J}$ is two.  Thus, $a_{i}$ appeared in both $I$ and $J$.  Since $a_{i+2}\in J$, this contradicts the assumption that $J$ is $r$-stable.  Similarly, in the latter case the multiplicity of $a_{i+2}$ in $\sort{I\cup J}$ is two, so $a_{i+2}$ must also appear in $I$, and this contradicts the assumption that $I$ is $r$-stable.    
Now suppose that $0<j<t$.  Then since $I$ is $r$-stable and $a_i\in I$, it must be that $a_{i+1}\in J$.  
But since $J$ is $r$-stable and $a_{i+2}\in J$, then $a_{i+1}\in I$, a contradiction.  
\qed
\bigskip

We let $\nabla_{n,k}^r$ denote the triangulation of $\rstablenk$ induced by $\nabla_{n,k}$.  This gives the following nesting of triangulations
\begin{equation*}
\nabla_{n,k}\supset\nabla_{n,k}^2\supset\nabla_{n,k}^3\supset\cdots\supset\nabla_{n,k}^{\left\lfloor\frac{n}{k}\right\rfloor},
\end{equation*}

\noindent thereby proving Theorem \ref{inducedtriangulationthm}.
In Section \ref{shelling}, the following lemma will play a key role.

\begin{lem}\label{fulldim}
If $n\equiv 1\mod k$, then $\khsr$ is $(n-1)$-dimensional for all $r\in\left[\lfloor\frac{n}{k}\rfloor\right]$.  In particular, $\Delta_{n,k}^{stab\left(\lfloor\frac{n}{k}\rfloor\right)}$ is a unimodular $(n-1)$-simplex.  
\end{lem}

\begin{proof}
Notice first that for $r=\left\lfloor\frac{n}{k}\right\rfloor$ there are precisely $n$ $r$-stable $k$-subsets of $[n]$.  Hence, $\rstablenk$ is an $(n-1)$-dimensional simplex.  Now suppose  $\epsilon$ is a vertex of this simplex.  Then precisely $k$ entries in $\epsilon$ are occupied by 1's, $k-1$ pairs of these 1's are separated by $r-1$ 0's, and the remaining pair is separated by $r$ 0's.  Hence, the only 1 that can be moved to the right and result in another $r$-stable vertex is the left-most 1 in the pair of 1's separated by $r$ 0's.  Making this move $n$ times results in returning to the vertex $\epsilon$, and produces a minimal circuit $(\omega)$ in $G_{n,k}$ using only $r$-stable vertices.  Since there are only $n$ such vertices it must be that $\sigma_{(\omega)}=\rstablenk$.  

We may also prove this result using Sturmfels' construction of this triangulation.  Simply notice that there are precisely $n$ $\left\lfloor\frac{n}{k}\right\rfloor$-stable $k$-subsets of $[n]$, namely

$$\begin{tabular}{c}
$\{1,1+\left\lfloor\frac{n}{k}\right\rfloor,1+2\left\lfloor\frac{n}{k}\right\rfloor,\ldots,1+(k-1)\left\lfloor\frac{n}{k}\right\rfloor\}$,\\
$\{2,2+\left\lfloor\frac{n}{k}\right\rfloor,2+2\left\lfloor\frac{n}{k}\right\rfloor,\ldots,2+(k-1)\left\lfloor\frac{n}{k}\right\rfloor\}$,\\
\vdots\\
$\{n,n+\left\lfloor\frac{n}{k}\right\rfloor,n+2\left\lfloor\frac{n}{k}\right\rfloor,\ldots,n+(k-1)\left\lfloor\frac{n}{k}\right\rfloor\}$.\\
\end{tabular}$$

It is easy to see that these subsets form a sorted collection of $k$-subsets of $[n]$.  Hence, they correspond to a unimodular $(n-1)$-simplex in the triangulation $\nabla_{n,k}$.
\end{proof}

In the coming sections we utilize the triangulation $\nabla_{n,k}^r$ to investigate geometric properties of the subpolytope $\rstablenk$ and its relationship with $\hypernk$.

%%%%%%%%%%%%%%%%%%%%%%%%%%%%%%%%%%%%%%%%%%%%%%%%%%%%%%%%%%%%%%%%%%%%%%%%%%%%%%%%%%
%%%%%%%%%%%%%%%%%%%%%%%%%%%%%%%%%%%%%%%%%%%%%%%%%%%%%%%%%%%%%%%%%%%%%%%%%%%%%%%%%%
%%%%%%%%%%%%%%%%%%%%%%%%%%%%%%%%%%%%%%%%%%%%%%%%%%%%%%%%%%%%%%%%%%%%%%%%%%%%%%%%%%

\section{A Shelling of the $r$-stable Second Hypersimplex}
\label{shelling}
Triangulations have many useful properties and well-studied applications in Ehrhart Theory.
Given a triangulation $\nabla$ of a $d$-dimensional polytope $\p$ let $\max\nabla$ denote the set of $d$-dimensional simplicies in $\nabla$.  
We call an ordering of the simplicies in $\max\nabla$, $(\alpha_1,\ldots,\alpha_s)$, a \emph{shelling} of $\nabla$ if for each $2\leq i\leq s$, $\alpha_i\cap(\alpha_1\cup\cdots\cup\alpha_{i-1})$ is a union of facets ($(d-1)$-dimensional faces) of $\alpha_i$.  
An equivalent condition for a shelling is that every $\alpha_i$ has a unique minimal (with respect to dimension) face that is not a face of the previous simplicies \cite{stan}.  
A triangulation with a shelling is called \emph{shellable}.  
For a shelling and a maximal simplex $\alpha$ in the triangulation define the \emph{shelling number} of $\alpha$, denoted $\#(\alpha)$, to be the number of facets shared by $\alpha$ and some previous simplex.  
The following theorem is due to Stanley.

\begin{thm}
\label{shellingthm}
\cite[Stanley]{stan3}
Let $\nabla$ be a unimodular shellable triangulation of a $d$-dimensional polytope $\p$.  Then 
\begin{equation*}
\sum_{j=0}^dh_j^*z^j=\sum_{\alpha\in\max\nabla}z^{\#(\alpha)}.
\end{equation*}
\end{thm}

In this section, we define a shelling of the triangulation $\nabla_{n,2}$ of the second hypersimplex $\hyperntwo$ that first shells the simplices within $\nabla_{n,2}^{r}$ and then extends this to a shelling of the $\nabla_{n,2}^{r-1}$ for every $0<r<\left\lfloor\frac{n}{2}\right\rfloor$, thereby proving Theorem~\ref{introshellingthm}.
%\begin{thm}\label{introshellingthm}
%There exists a shelling of $\nabla_{n,2}$ that first builds $\Delta_{n,2}^{stab(r)}$ and then builds $\Delta_{n,2}^{stab(r-1)}$ for every $0<r\leq\left\lfloor\frac{n}{2}\right\rfloor$.  
%%Hence, we say there exists a {\bf stable shelling} of the odd second hypersimplex.  
%\end{thm}
%\noindent Theorem \ref{introshellingthm} follows from Theorem \ref{introshellingthm}.
The following remark outlines our proof of Theorem~\ref{introshellingthm}.

\begin{rmk}
\label{proofapproach}
We will prove Theorem~\ref{introshellingthm} in two cases, when $n$ is odd and when $n$ is even.  
The bulk of the work will be done in the case when $n$ is odd, and then we will quickly extend to the even case.  

Fix $k=2$ and $n$ odd.  
Notice first that by Lemma \ref{fulldim} we can certainly shell $\nabla_{n,2}^r$ where $r=\left\lfloor\frac{n}{2}\right\rfloor$.  
Our goal is to inductively shell $\nabla_{n,2}^r$ with this as our base case.
That is, assume we have previously shelled $\nabla_{n,2}^{r+1}$ for $0<r<\left\lfloor\frac{n}{2}\right\rfloor$.  
We then describe a continuation of this shelling to $\nabla_{n,2}^r$.
Since we are assuming we have previously shelled the simplices in $\nabla_{n,2}^{r+1}$ we must extend this order to the set of simplices $\omega\in\max\nabla_{n,2}^r\backslash\max\nabla_{n,2}^{r+1}$.  Each simplex $\omega$ in this set uses some vertices that are $r$-stable but not $(r+1)$-stable.  We will call these vertices the \emph{$r$-adjacent vertices}.  We first select a particular $r$-adjacent vertex of $\omega$, and think of this as the initial vertex in the cycle $(\omega)$.  Given this choice, we then associate to $\omega$ a composition of $r$ into $n-r-1$ parts that describes the cycle $(\omega)$ in terms of the selected initial vertex.  Using this composition and its relationship with the vertices of $\omega$ we associate to $\omega$ a lattice path in a decorated ladder-shaped region of the plane.  We then order the simplices in $\omega\in\max\nabla_{n,2}^r\backslash\max\nabla_{n,2}^{r+1}$ by first collecting them into sets based on their initial vertex and the number of $r$-adjacent vertices they use, ordering these sets from least to most $r$-adjacent vertices used, and then ordering the elements within these sets via the colexicographic ordering applied to their associated compositions.  We then utilize their associated lattice paths to identify the unique minimal new face for each simplex.  In particular, we shell the simplices in terms of least $r$-adjacent vertices used to most $r$-adjacent vertices used.

Finally, when $n$ is even, we simply adjust the base case and then the rest of the results will extend naturally.  
For $n$ even, the base case will be a shelling of $\nabla_{n,2}^r$ for $n=2r+2$ and $r\geq1$.  
\end{rmk}

Following Remark~\ref{proofapproach}, we now fix $k=2$ and $n$ odd, until we reach subsection~\ref{the n even case}.  

\subsection{$r$-adjacent vertices}

For $\ell\in[n]$, let $\adj{\ell}$ denote the vertex $\epsilon_I$ where $I=\{\ell,\ell+r\}\in{[n]\choose 2}$.  We call a vertex $\adj{\ell}$ an \emph{$r$-adjacent vertex}.
Let $\Adj{n}:=\left\{\adj{\ell}:\ell\in[n]\right\}$.  So $\Adj{n}$ is precisely the set of vertices that are $r$-stable but not $(r+1)$-stable.

\begin{lem}\label{reachablevertices}
Let $\epsilon$ and $\epsilon^\prime$ be two vertices in $(\omega)$ for some simplex $\omega\in\max\nabla_{n,2}$.  Suppose that $\epsilon$ has entries $\epsilon_i=\epsilon_j=1$ with $i<j$, and $\epsilon_t=0$ for all $t\neq i,j$.  Suppose also that $\epsilon^\prime$ has entries $\epsilon^\prime_k=\epsilon^\prime_l=1$ with $k<l$, and $\epsilon^\prime_t=0$ for all $t\neq k,l$.  Then (modulo $n$) we have that
$$i\leq k\leq j \leq l.$$
\end{lem}

\begin{proof}
Since $\epsilon$ and $\epsilon^\prime$ are vertices of a simplex in $\nabla_{n,2}$ they correspond to a sorted pair of $2$-subsets of $[n]$.
\end{proof}

\begin{cor}\label{circdistadj}
Let $\omega\in\max\nabla_{n,2}$, and suppose that $\adj{\ell}$ and $\adj{\ell^\prime}$ are vertices in $v_{(\omega)}\cap\Adj{n}$.  Then $\cd{\ell}{\ell^\prime}\leq r$.  
\end{cor}

\begin{proof}
Lemma \ref{reachablevertices} indicates that 
\begin{equation*}
\begin{split}
\ell\leq\ell^\prime&\leq\ell+r\leq\ell^\prime+r, \mbox{ or}\\
\ell^\prime\leq\ell&\leq\ell^\prime+r\leq\ell+r.\\
\end{split}
\end{equation*}
\end{proof}

\begin{rmk}\label{orderadj}
Consider a simplex $\omega\in\max\nabla_{n,2}$.  Notice that for a fixed $0<r< \left\lfloor\frac{n}{2}\right\rfloor$ we may order the elements of the set $v_{(\omega)}\cap\Adj{n}$ as $\adj{\ell}<_{adj}\adj{\ell^\prime}$ if and only if $\ell<\ell^\prime$.  In this way, there exists a unique maximal element of the set $v_{(\omega)}\cap\Adj{n}$.
\end{rmk}

\subsection{Associating a composition to $\omega\in\max\nabla_{n,2}^r\backslash\max\nabla_{n,2}^{r+1}$}
Fix $\omega\in\max\nabla_{n,2}^r\backslash\max\nabla_{n,2}^{r+1}$.  Then $\omega$ uses at least one element of $\Adj{n}$.  We may fix one such $\adj{\ell}$, and consider the circuit $(\omega)$ as having initial vertex $\adj{\ell}$.  We refer to the $1$ in entry $\ell$ of the vertex $\adj{\ell}$ as the \emph{left $1$} and the $1$ in entry $\ell+r$ as the \emph{right $1$}.  Then, in $(\omega)$ each edge corresponds to a move of the left 1 or of the right 1.  In particular, 
\begin{enumerate}[($\star$)]
	\item the left 1 makes $r$ moves,
	\item the right 1 makes $n-r$ moves, and
	\item the left 1 cannot move first or last.
\end{enumerate}
Note that the first two conditions are immediate from the definition of $(\omega)$ and the fact that $k=2$.  The third condition holds since $\omega$ uses only vertices that are $r$-stable.  It follows that for a fixed $\adj{\ell}\in v_{(\omega)}$ we can think of the circuit $(\omega)$ as a sequence of moves of the left 1 and moves of the right 1 satisfying these conditions.  Moreover, we may encode this as a composition
			$$\lambda=(\lambda_1,\lambda_2,\ldots,\lambda_{n-r-1})$$
of $r$ into $n-r-1$ parts, where part $\lambda_i$ denotes the number of moves of the left 1 after the $i^{th}$ move of the right 1 and before the $(i+1)^{st}$ move of the right 1.

\begin{example}\label{compositionexample(15)4}
Consider the minimal circuit in $G_{15,2}$ depicted in Figure \ref{fig:minimalcircuitexample(15)4}.  This circuit corresponds to a simplex $\omega\in\max\nabla_{15,2}^4\backslash\max\nabla_{15,2}^{5}$.  If we choose the initial vertex of this circuit to be the unique maximal element of the set $v_{(\omega)}\cap\Adjvar{4}{15}$, namely $\adjvar{4}{15}$, then this circuit has associated composition $\lambda=(1,0,0,1,0,1,0,0,0,1).$

\begin{figure}
	\centering
	\includegraphics[width=0.9\textwidth]{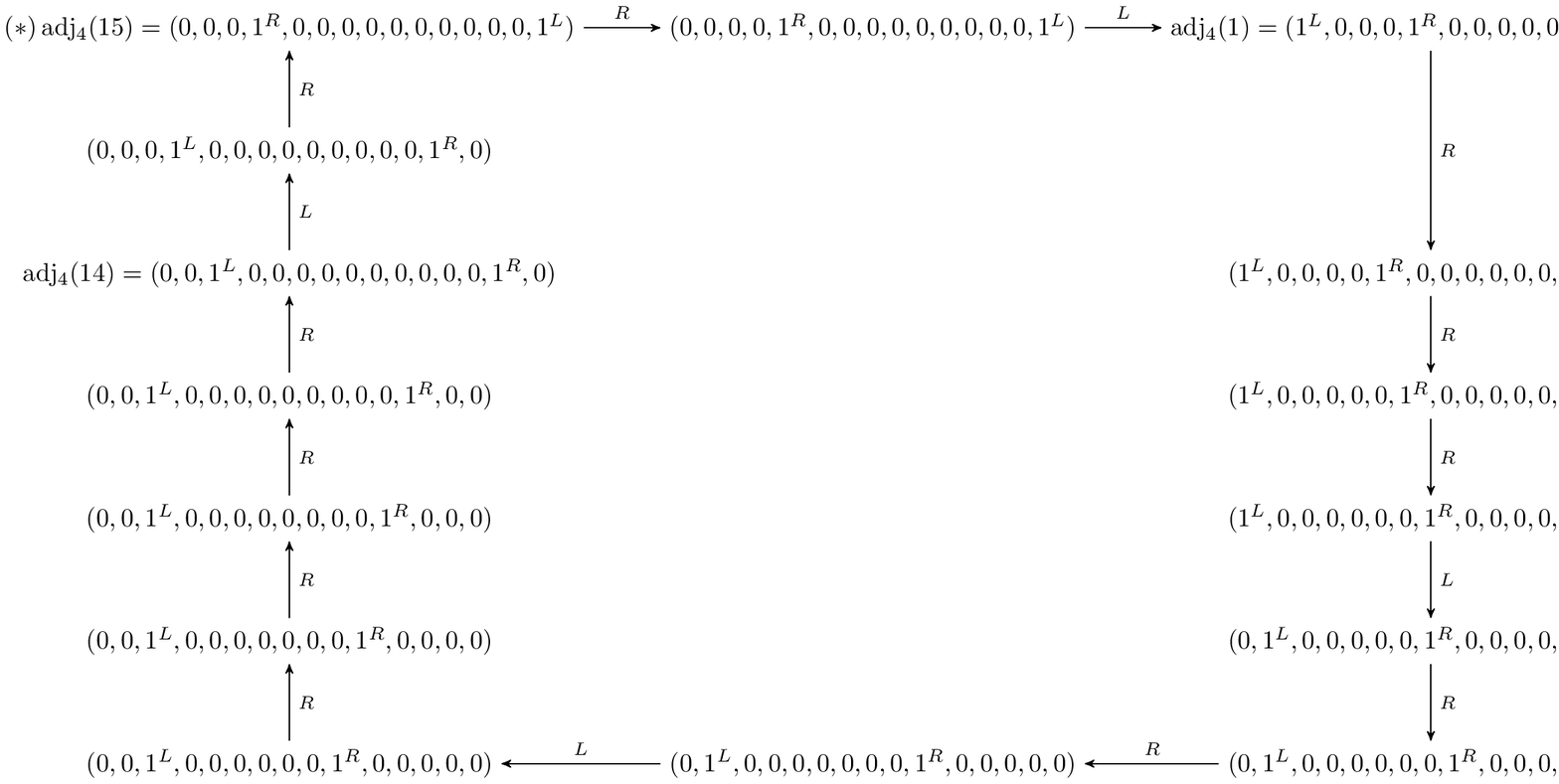}
	\caption{The minimal circuit corresponding to the simplex 
	$$\omega=5671892(10)(11)(12)(13)3(14)4(15).$$}
	\label{fig:minimalcircuitexample(15)4}
\end{figure}

\end{example}

\begin{example}\label{compositionexample93}
Next consider the minimal circuit in $G_{9,2}$ depicted in Figure \ref{fig:minimalcircuitexample93}.  This circuit corresponds to a simplex $\omega\in\max\nabla_{9,2}^3\backslash\max\nabla_{9,2}^{4}$.  If we choose the initial vertex of this circuit to be the unique maximal element of the set $v_{(\omega)}\cap\Adjvar{3}{9}$, namely $\adjvar{3}{7}$, then this circuit has associated composition $\lambda=(0,0,3,0,0).$

\begin{figure}
	\centering
	\includegraphics[width=0.9\textwidth]{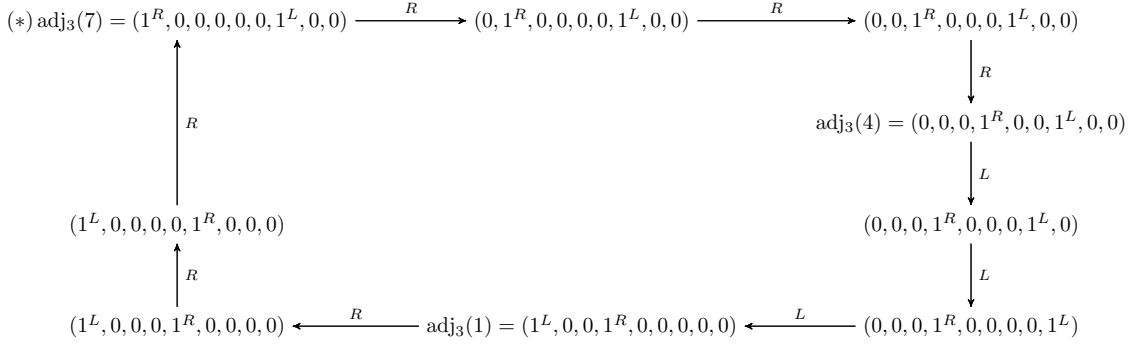}
	\caption{The minimal circuit corresponding to the simplex 
	$\omega=456123789.$}
	\label{fig:minimalcircuitexample93}
\end{figure}

\end{example}

\begin{prop}\label{compositioncorrespondence}
Fix $\adj{\ell}\in\Adj{n}$.  Each simplex $\omega\in \max\nabla_{n,2}^r\backslash\max\nabla_{n,2}^{r+1}$ that uses the vertex $\adj{\ell}$ corresponds uniquely to a composition
			$$\lambda=(\lambda_1,\lambda_2,\ldots,\lambda_{n-r-1})$$
of $r$ into $n-r-1$ parts that satisfies
\begin{equation}\label{lambdabounds}
i+1+2r-n\leq\sum_{j=1}^i\lambda_j\leq i
\end{equation}
for all $i=1,2,\ldots,n-r-1$.
\end{prop}

\begin{proof}
Let $\omega\in \max\nabla_{n,2}^r\backslash\max\nabla_{n,2}^{r+1}$ that uses vertex $\adj{\ell}$.  Then $(\omega)$ is a minimal circuit in the directed graph $G_{n,2}$, one of whose vertices is $\adj{\ell}$.  Thinking of $\adj{\ell}$ as the initial vertex we consider the 1 in place $\ell$ as the left 1 and the 1 in place $\ell+r$ as the right 1.  By the above conditions it is clear that we may construct the partition $\lambda$ of  $r$ into $n-r-1$ parts, where part $\lambda_i$ denotes the number of moves of the left 1 after the $i^{th}$ move of the right 1 and before the $(i+1)^{st}$ move of the right 1.  Since $\omega\in\max\nabla_{n,2}^r$ the left 1 can never have made more moves that the right 1.  This gives the upper bound on $\sum_{j=1}^i\lambda_j$ for each $i=1,2,\ldots,n-r-1$.

Similarly, since $\omega\in\max\nabla_{n,2}^r$ after the $(r+1)^{st}$-to-last move of the right 1 and before its $r^{th}$-to-last move we must have that the left 1 moved at least once.  More generally, after the $n-r-t+1^{st}$ move of the right 1 we must have that the left 1 moved at least $r-t+2$ times for $t=r+1,r,r-1,\ldots,2$.  Hence, the number of left moves that occur after the $i^{th}$ right move and before the $(i+1)^{st}$ right move is at least $i+1+2r-n$.  This gives the lower bound on $\sum_{j=1}^i\lambda_j$.  

Conversely, suppose that we have a composition $\lambda$ satisfying the given conditions.  We can construct a simplex $\omega(\lambda)\in\max\nabla_{n,2}^r\backslash\max\nabla_{n,2}^{r+1}$ that uses the vertex $\adj{\ell}$ by constructing a minimal circuit in $G_{n,2}$ as follows.  Starting with $\adj{\ell}$, and labeling the left 1 and right 1 as we have been, after the $i^{th}$ move of the right 1 move the left 1 $\lambda_i$ times.  Once this has been done for all $i=1,2,\ldots,n-r-1$, move the right 1 once more.  The upper bound ensures that the right distance between the 1s is always at least $r$.  Similarly, the lower bound ensures that the left distance is always at least $r$.  Since $\adj{\ell}$ is in the circuit $(\omega(\lambda))$ then this corresponds to a simplex $\omega(\lambda)\in\max\nabla_{n,2}^r\backslash\max\nabla_{n,2}^{r+1}$.
\end{proof}

\begin{rmk}
\label{compositionlabel}
By Remark \ref{orderadj} we can identify each simplex $\omega\in\max\nabla_{n,2}^r\backslash\max\nabla_{n,2}^{r+1}$ with the unique maximal element of $v_{(\omega)}\cap\Adj{n}$, say $\adj{\ell}$.  Let $\lambda$ be the composition associated to $\omega$ via $\adj{\ell}$ by Proposition \ref{compositioncorrespondence}.  Then we may uniquely label the simplex $\omega$ as $\omega_{\ell,\lambda}$.
\end{rmk}

\begin{defn}
\label{parity}
For a simplex $\omega_{\ell,\lambda}$ recall that we think of the 1 in entry $\ell$ of $\adj{\ell}$ as the left 1, and the 1 in entry $\ell+r$ as the right 1.  
	\begin{itemize}
		\item A {\bf left move} in $(\omega_{\ell,\lambda})$ is an edge in $(\omega_{\ell,\lambda})$ corresponding to a move of the left 1, and
		\item A {\bf right move} in $(\omega_{\ell,\lambda})$ is an edge in $(\omega_{\ell,\lambda})$ corresponding to a move of the right 1.  
		\item The {\bf parity} of an edge in $(\omega_{\ell,\lambda})$ is \emph{left} if the edge is a left move, and \emph{right} if the edge is a right move.
	\end{itemize}
\end{defn}

\begin{rmk}[{\bf Lattice Path Correspondence}]
\label{latticepath}
Notice that each simplex $\omega_{\ell,\lambda}$ in the set $\max\nabla_{n,2}^r\backslash\max\nabla_{n,2}^{r+1}$ corresponds to a lattice path, $p(\omega_{\ell,\lambda})$, from $(0,0)$ to $(n-r,r)$ that uses only North (0,1) and East (1,0) moves.  Here, right moves in the circuit $(\omega_{\ell,\lambda})$ correspond to East moves in $p(\omega_{\ell,\lambda})$, and left moves in  $(\omega_{\ell,\lambda})$ correspond to North moves in $p(\omega_{\ell,\lambda})$.  By Proposition \ref{compositioncorrespondence} the lattice path $p(\omega_{\ell,\lambda})$ is bounded between the lines $y=x$ and $y=x-n+2r$.  Each vertex in $(\omega_{\ell,\lambda})$ corresponds uniquely to a lattice point on $p(\omega_{\ell,\lambda})$.  In particular, for $0\leq t\leq r$ the vertex $\adj{\ell+t}$ corresponds to the lattice point $(t,t)$, and the vertex $\adj{\ell-r+t}$ corresponds to the lattice point $(t,n-2r+t)$.  
\end{rmk}

Here are some examples of simplices and their corresponding lattice paths.

\begin{example}\label{ex:latticepathexample(15)4}
Let $n=15$ and $r=4$.  Recall the simplex from Example \ref{compositionexample(15)4}
	$$\omega=5671892(10)(11)(12)(13)3(14)4(15)\in\max\nabla_{15,2}^4\backslash\max\nabla_{15,2}^{5}.$$
This simplex corresponds to the minimal circuit in the graph $G_{15,2}$ depicted in Figure \ref{fig:minimalcircuitexample(15)4}.

From this, we can see that $\omega$ uses the vertices $\adjvar{4}{15}$, $\adjvar{4}{1}$, and $\adjvar{4}{14}$.  Hence, we label $\omega$ as $\omega_{15,\lambda}$, where
	$$\lambda=(1,0,0,1,0,1,0,0,0,1).$$
The lattice path corresponding to $\omega$ via this labeling is depicted in Figure \ref{fig:latticepathexample(15)4}.

	\begin{figure}
		\centering
		\includegraphics[width=0.7\textwidth]{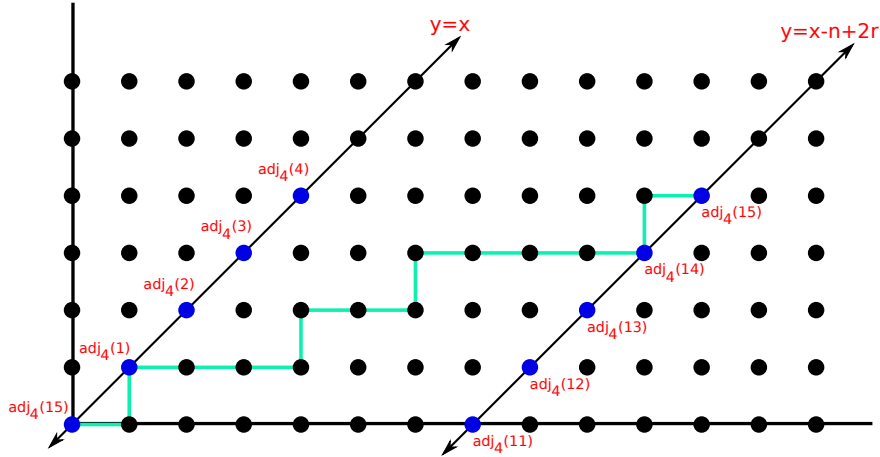}
		\caption{The lattice path $p(\omega_{15,\lambda})$, where $\lambda=(1,0,0,1,0,1,0,0,0,1)$.}
		\label{fig:latticepathexample(15)4}
	\end{figure}
\end{example}

\begin{example}\label{ex:latticepathexample93}
Let $n=9$ and $r=3$.  Recall the simplex from Example \ref{compositionexample93}
	$$\omega=456123789\in\max\nabla_{9,2}^3\backslash\max\nabla_{9,2}^{4}.$$
This simplex corresponds to the minimal circuit in the graph $G_{9,2}$ depicted in Figure \ref{fig:minimalcircuitexample93}.

From this, we can see that $\omega$ uses the vertices $\adjvar{3}{7}$, $\adjvar{3}{4}$, and $\adjvar{3}{1}$.  Hence, we label $\omega$ as $\omega_{7,\lambda}$, where
	$$\lambda=(0,0,3,0,0).$$
The lattice path corresponding to $\omega$ via this labeling is depicted in Figure \ref{fig:latticepathexample93}.

	\begin{figure}
		\centering
		\includegraphics[width=0.5\textwidth]{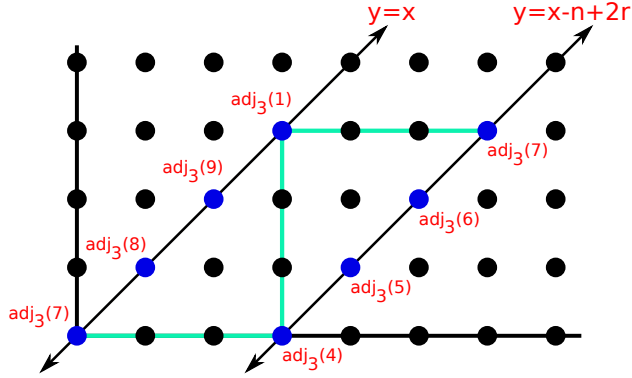}
		\caption{The lattice path $p(\omega_{7,\lambda})$, where $\lambda=(0,0,3,0,0)$.}
		\label{fig:latticepathexample93}
	\end{figure}
\end{example}

\subsection{The shelling order}
\label{the shelling order}
Recall that the \emph{colexicographic} order on a pair of ordered $m$-tuples ${\bf a}=(a_1,\ldots,a_m)$ and ${\bf b}=(b_1,\ldots,b_m)$ is defined by ${\bf b}<_{\colex}{\bf a}$ if and only if the right-most nonzero entry in ${\bf a}-{\bf b}$ is positive.  Let $W_{\ell,s}$ denote the set of all simplices with label $\omega_{\ell,\lambda}$ that use precisely $s$ elements of $\Adj{n}$.  
  Order the elements in each set $W_{\ell,s}$ with respect to the colexicographic ordering on their associated compositions (from least to greatest).  We write $\omega_{\ell,\lambda}<_{\colex}\omega_{\ell,\lambda^\prime}$ if and only if $\lambda<_{\colex}\lambda^\prime$.  Next order the sets $W_{\ell,s}$ (from least to greatest) with respect to the colexicographic ordering on the labels $(\ell,s)$.  We then write $\omega_{\ell,\lambda}<\omega_{\ell^\prime,\lambda^\prime}$ if and only $\omega_{\ell,\lambda}\in W_{\ell,s}$ and $\omega_{\ell^\prime,\lambda^\prime}\in W_{\ell^\prime,s^\prime}$ with $(\ell,s)<_{\colex}(\ell^\prime,s^\prime)$ or if $(\ell,s)=(\ell^\prime,s^\prime)$ and $\omega_{\ell,\lambda}<_{\colex}\omega_{\ell,\lambda^\prime}$.

\begin{thm}\label{stableshellingthm}
The order $<$ on the simplices $\omega\in\max\nabla_{n,2}^r\backslash\max\nabla_{n,2}^{r+1}$ (from least to greatest) extends the shelling of $\nabla_{n,2}^{r+1}$ to a shelling of $\nabla_{n,2}^r$.  
\end{thm}

The $n$ odd case of Theorem \ref{introshellingthm} follows immediately from Theorem \ref{stableshellingthm}.  
To prove Theorem \ref{stableshellingthm} it suffices to identify the unique minimal new face associated to each simplex in the shelling order.  
To do so, we first prove a sequence of lemmas.

\begin{lem}\label{adjminnotmaxmoves}
Suppose the $\omega_{\ell,\lambda}$ uses $\adj{\ell^\prime}$ for $\ell^\prime\neq\ell$.  Then $\adj{\ell^\prime}$ is a vertex in $(\omega_{\ell,\lambda})$ that is either
\begin{enumerate}[(i)]
	\item produced by a right move for which the preceding number of left moves is minimal and not maximal with respect to equation (\ref{lambdabounds}), or 
	\item produced by a left move and followed by a right move for which the number of left moves preceding the right move is maximal and not minimal with respect to equation (\ref{lambdabounds}).
\end{enumerate}
\end{lem}

\begin{proof}
Since $\adj{\ell}$ is selected to be the greatest element of $v_{(\omega)}\cap\Adj{n}$ then each other $\adj{\ell^\prime}$ used by $\omega_{\ell,\lambda}$ is produced in $(\omega_{\ell,\lambda})$ by doing $n-2r+t$ right moves for some number $t$ of left moves, or $\adj{\ell^\prime}$ is produced by doing $0<t<r$ right moves and the same number of left moves.  In the former case, such a vertex corresponds to an entry $\lambda_m$ in the composition $\lambda$ with $m=n-2r+t$ for which
	$$t=m+2r-n=(m-1)+1+2r-n\leq\sum_{j=1}^{m-1}\lambda_j=t.$$
Hence, the number of left moves preceding the $m^{th}$ right move is minimal.  Moreover, the number of left moves preceding the $m^{th}$ right move is maximal only if
\begin{equation*}
m+2r-n=t=\sum_{j=1}^{m-1}\lambda_j=m-1.
\end{equation*}
Thus,
\begin{equation*}
r=\frac{n-1}{2}=\left\lfloor\frac{n}{2}\right\rfloor.
\end{equation*}
But recall that since $\Delta_{n,2}^{stab\left(\left\lfloor\frac{n}{2}\right\rfloor\right)}$ is a unimodular $(n-1)$-simplex we are only completing the shelling of $\nabla_{n,2}^{r+1}$ to a shelling of $\nabla_{n,2}^r$ for $r<\left\lfloor\frac{n}{2}\right\rfloor$.  So we conclude that the sum is minimal and not maximal.

In the latter case, the vertex $\adj{\ell_j}$ is produced by doing $0<t<r$ right moves and the same number of left moves.  Thus, following this vertex with another left move results in a vertex that is no longer $r$-stable.  So the move following $\adj{\ell_j}$ must be a right move.  Such a vertex corresponds to an entry $\lambda_m$ in the composition $\lambda$ for which $m=t$, and the right move following the vertex is the $(m+1)^{st}$ right move in the circuit.  Thus,
$$m+1+2r-n\leq\sum_{j=1}^m\lambda_j=t=m.$$
Hence, the number of left moves preceding the right move following the vertex is maximal.  If this number is also minimal then it must be that 
\begin{equation*}
\begin{split}
m+1+2r-n&=m,\\
r&=\frac{n-1}{2}=\left\lfloor\frac{n}{2}\right\rfloor,\\
\end{split}
\end{equation*}
and so we conclude that the sum is not minimal just as in the previous case.  It remains to show that the move preceding the vertex $\adj{\ell_j}$ is a left move.  Suppose for the sake of contradiction that $\adj{\ell_j}$ is preceded by a right move.  Then $\lambda_m=0$.  Thus, since the number of left moves preceding the $(m+1)^{st}$ right move is maximal we have that
$$m=\sum_{j=1}^m\lambda_j=\sum_{j=1}^{m-1}\lambda_j\leq m-1,$$
which is a contradiction.  Thus, we conclude that $\adj{\ell_j}$ is produced by a left move and followed by a right move for which the number of left moves preceding the right move is not minimal.
\end{proof}

\begin{lem}\label{parityforadj}
Suppose that the simplex $\omega_{\ell,\lambda}$ uses the elements 
	$$\adj{\ell_1}<_{adj}\adj{\ell_2}<_{adj}\cdots<_{adj}\adj{\ell_s}=\adj{\ell}$$
of $\Adj{n}$.  For $j\neq s$, the parities of the edges preceding $\adj{\ell_j}$ in $(\omega_{\ell,\lambda})$ and following $\adj{\ell_j}$ are opposite.  Also, the parity of the edges about $\adj{\ell}$ is right.
\end{lem}

\begin{proof}
First recall that we have already noted that the first and last moves of $(\omega_{\ell,\lambda})$ must be right moves.  Hence, the parity of the edges about $\adj{\ell}$ is right.  

Now consider $\adj{\ell_j}$ for $j\neq s$.  By Lemma \ref{adjminnotmaxmoves} we have two cases.  In case (ii), $\adj{\ell_j}$ is produced by a left move and followed by a right move for which the number of left moves preceding the right move is not minimal.  Hence, the result is immediate.  

In case (i), $\adj{\ell_j}$ is produced by a right move for which the preceding number of left moves is minimal and not maximal.  Suppose for the sake of contradiction that the parities of the moves about $\adj{\ell_j}$ are the same.  So if $\adj{\ell_j}$ is produced by the $m^{th}$ right move we have that
	$$m+2r-n=(m-1)+1+2r-n=\sum_{j=1}^{m-1}\lambda_j.$$
Since the parities of the edges about $\adj{\ell_j}$ are the same it is followed by a right move, and so it must be that $\lambda_m=0$.  Hence, by equation (\ref{lambdabounds})
	$$m+1+2r-n\leq\sum_{j=1}^m\lambda_j=\sum_{j=1}^{m-1}\lambda_j=m+2r-n,$$
which is a contradiction.
\end{proof}

\begin{lem}\label{switchingparityforadj}
Suppose that the simplex $\omega_{\ell,\lambda}$ uses the vertex $\adj{\ell^\prime}$.  Then switching the parities of the moves about $\adj{\ell^\prime}$ does not replace $\adj{\ell^\prime}$ with another vertex in $\Adj{n}$.  
\end{lem}

\begin{proof}
First consider the case where $\ell^\prime\neq\ell$.  By Remark \ref{latticepath} the simplex $\omega_{\ell,\lambda}$ corresponds to a lattice path $p(\omega_{\ell,\lambda})$ that is bounded between the lines $y=x$ and $y=x-n+2r$, and the elements of $\Adj{n}$ reachable from $\adj{\ell}$ all lie on these two lines.  Suppose for the sake of contradiction that switching the parities of the moves about $\adj{\ell^\prime}$ replaced this vertex with another element of $\Adj{n}$, say $\adj{\ell^{\prime\prime}}$.  Then the resulting change in the lattice path $p(\omega_{\ell,\lambda})$ implies that $\adj{\ell^{\prime\prime}}$ lies on the opposite of these two lines from that on which $\adj{\ell^\prime}$ lies.  It then follows that $n-2r=2$, or equivalently, $n=2r+2$.  Since we have chosen $n$ to be odd this is a contradiction.

Now consider the case where $\ell^\prime=\ell$.  Suppose for the sake of contradiction that switching the parities of the moves about $\adj{\ell}$ replaces $\adj{\ell}$ with another vertex $\adj{\ell^{\prime\prime}}$.  Consider the vertex before the right move producing $\adj{\ell}$.  Since this right move produces $\adj{\ell}$ then in this preceding vertex there must be precisely $r$ $0$'s to the right of the right 1 and before the left 1.  Similarly, since the left move produces the vertex $\adj{\ell^{\prime\prime}}$ it must be that there are $r$ $0$'s to the right of the left 1 and before the right 1.  Hence, $n=2r+2$, a contradiction.
\end{proof}

\begin{lem}\label{replacingadj}
Suppose that the simplex $\omega_{\ell,\lambda}$ uses the elements 
	$$\adj{\ell_1}<_{adj}\adj{\ell_2}<_{adj}\cdots<_{adj}\adj{\ell_s}=\adj{\ell}$$
of $\Adj{n}$.  Switching the parity of the edges about $\adj{\ell_j}$ replaces the vertex $\adj{\ell_j}$ with an $(r+1)$-stable vertex, and leaves all other vertices in $(\omega_{\ell,\lambda})$ fixed.
\end{lem}

\begin{proof}
First fix $\adj{\ell_j}$ for $j\neq s$, and switch the parity of the moves directly before and after $\adj{\ell_j}$ in $(\omega_{\ell,\lambda})$.  By Lemma \ref{adjminnotmaxmoves} there are two cases.  In case (i), Lemma \ref{parityforadj} implies that the parity switch changes the move before from a right to a left, and the move after from a left to a right.  Since each vertex in the circuit is determined by the number of left moves and right moves by which it differs from $\adj{\ell}$ this switching does not change any of the vertices preceding $\adj{\ell_j}$ in $(\omega_{\ell,\lambda})$.  Similarly, it does not change any of the vertices following $\adj{\ell_j}$.  The reader should also note that this switch changes the composition $\lambda$.  However, Lemma \ref{adjminnotmaxmoves} ensures that the resulting composition, say $\lambda^\prime$, still satisfies the bounds of equation (\ref{lambdabounds}).  Hence, by Proposition \ref{compositioncorrespondence} this switch produces a circuit $(\omega_{\ell,\lambda^\prime})$ for which $\omega_{\ell,\lambda^\prime}\in\max\nabla_{n,2}^r\backslash\max\nabla_{n,2}^{r+1}$.  Moreover, the vertex which replaces $\adj{\ell_j}$ is not an element of $\Adj{n}$ by Lemma \ref{switchingparityforadj}.  As an example, consider the following scenario for $r=3$:

\begin{center}
	\includegraphics[width=0.7\textwidth]{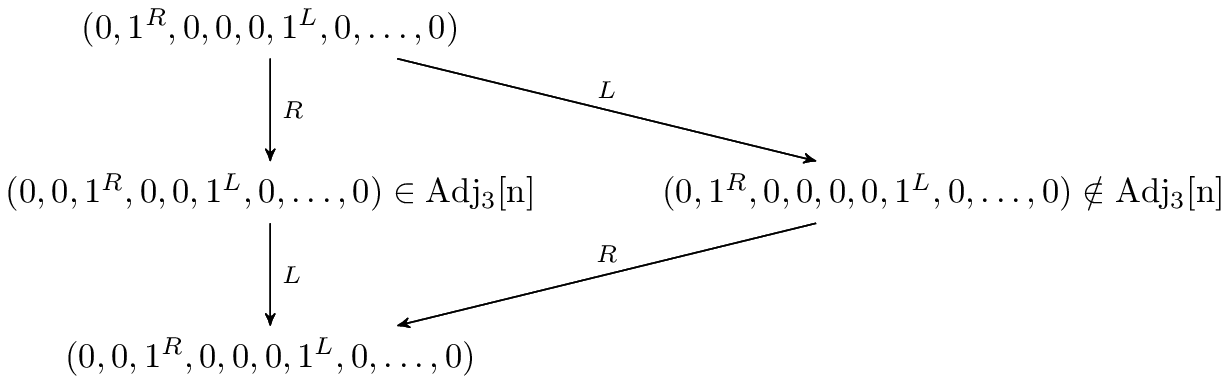}
\end{center}

We remark that $\omega_{\ell,\lambda}\in W_{\ell,s}$, so $\omega_{\ell,\lambda^\prime}\in W_{\ell,s-1}$.  

In case (ii), Lemma \ref{parityforadj} implies that the parity switch changes the move before $\adj{\ell_j}$ from a left to a right, and the move after $\adj{\ell_j}$ from a right to a left.  Now apply the same argument as for case (i), and the result follows.  We again remark that since $\omega_{\ell,\lambda}\in W_{\ell,s}$ then the parity switch results in a simplex $\omega_{\ell,\lambda^\prime}\in W_{\ell,s-1}$.

Now consider $\adj{\ell}$.  By the same argument as before, switching the parities of the moves about $\adj{\ell}$ replaces $\adj{\ell}$ with an $(r+1)$-stable vertex that is not in $\Adj{n}$.  This scenario is depicted in the following diagram for $r=3$.  

\begin{center}
	\includegraphics[width=0.6\textwidth]{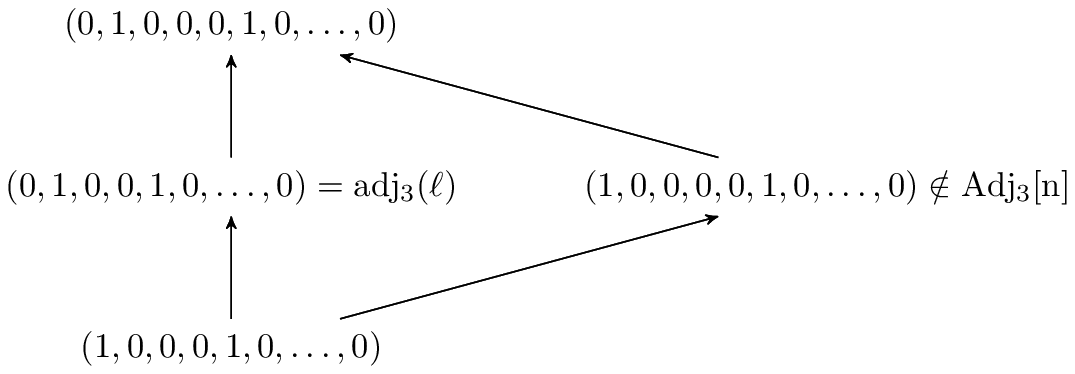}
\end{center}

Notice that we omit the labels $R$ and $L$.  This is because removing $\adj{\ell}$, the vertex which defines the labels, demands a relabeling of the resulting circuit, and in general the new labels will not agree with the old.  However, this is acceptable since to switch the parities of the moves about $\adj{\ell}$ we simply note that the vertices directly before and after $\adj{\ell}$ in the circuit are completely determined by $\adj{\ell}$.  Moreover, the edges before and after $\adj{\ell}$ each correspond to a move of a different 1 in $\adj{\ell}$.  Hence, to switch the parities, starting at the vertex preceding $\adj{\ell}$ simply switch the order in which the 1's move.

\end{proof}

\begin{cor}\label{cortoreplacingadj}
For the simplex $\omega_{\ell,\lambda}$, switching the parities of the edges about $\adj{\ell_j}$, for $j\neq s$, reduces the simplex $\omega_{\ell,\lambda}\in W_{\ell,s}$ to a simplex $\omega_{\ell,\lambda^\prime}\in W_{\ell,s-1}$.   For $s\neq 1$, switching the parities of the edges about $\adj{\ell}$ reduces the simplex $\omega_{\ell,\lambda}\in W_{\ell,s}$ to a simplex $\omega_{\ell_{s-1},\lambda^\prime}\in W_{\ell_{s-1},s-1}$.  For $s=1$, switching the parities of the edges about $\adj{\ell}$ reduces the simplex $\omega_{\ell,\lambda}\in W_{\ell,s}$ to a simplex in $\max\nabla_{n,2}^{r+1}$.
\end{cor}

\subsubsection{Proof of Theorem \ref{stableshellingthm} for $n$ odd}  
\label{stable shelling in odd case proof}
\begin{figure}
	\centering
	\includegraphics[width=0.7\textwidth]{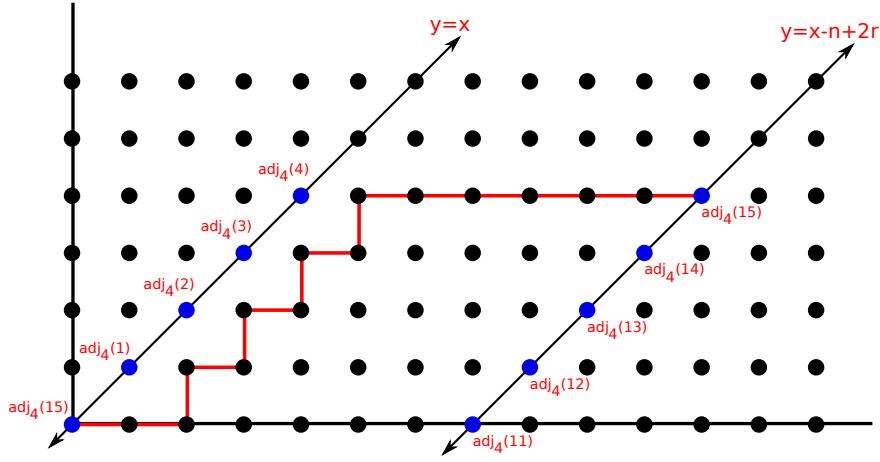}
	\caption{The lattice path $p(\omega_{15,\lambda^\star})$.  Here, $\lambda^\star=(0,1,1,1,1,0,0,0,0,0)$.}
	\label{fig:lambdastarexample(15)4}
\end{figure}

We are now ready to prove Theorem \ref{stableshellingthm} when $n$ is odd.  
Recall, to prove Theorem \ref{stableshellingthm} it suffices to identify the unique minimal new face associated to each simplex in the shelling order.  
Given $\omega_{\ell,\lambda}\in W_{\ell,s}$ recall that we can associate to $\omega_{\ell,\lambda}$ a lattice path $p(\omega_{\ell,\lambda})$.  
Notice also that for each $\adj{\ell}\in\Adj{n}$ the first simplex in our order that uses $\adj{\ell}$ is $\omega_{\ell,\lambda^{\star}}$ where
	$$\lambda^\star=(0,1,1,\ldots,1,0,0,\ldots,0).$$
A picture of the lattice path $p(\omega_{15,\lambda^\star})$ corresponding to $\lambda^\star$ for $n=15$, $r=4$, and $\ell=15$ is given in Figure \ref{fig:lambdastarexample(15)4}.

We claim that the unique minimal new face of $\omega_{\ell,\lambda}$ is the collection of vertices

\begin{itemize}
	\item $\adj{\ell}$,
		
	\item those vertices corresponding to lattice points on $p(\omega_{\ell,\lambda})$ that lie on $y=x$, and
		
	\item those vertices corresponding to lattice points on $p(\omega_{\ell,\lambda})$ that are
		\begin{itemize}
			\item right-most in their row of the lattice,
				
			\item corners of $p(\omega_{\ell,\lambda})$, and
				
			\item do not lie on $p(\omega_{\ell,\lambda^\star})$.
		\end{itemize}
\end{itemize}

That is to say, the corners of $p(\omega_{\ell,\lambda})$ that are ``furthest away" or ``point away" from the path $p(\omega_{\ell,\lambda^\star})$.

\begin{example}
Let $n=15$ and $r=4$.  Consider the simplex $\omega_{15,\lambda}$ from Example \ref{ex:latticepathexample(15)4}.  The unique minimal new face for $\omega_{15,\lambda}$ is given by the vertices
$$	\begin{array}{c}
(0,0,0,1,0,0,0,0,0,0,0,0,0,0,1)\\
(1,0,0,0,1,0,0,0,0,0,0,0,0,0,0)\\
(1,0,0,0,0,0,0,1,0,0,0,0,0,0,0)\\
(0,1,0,0,0,0,0,0,0,1,0,0,0,0,0)\\
(0,0,1,0,0,0,0,0,0,0,0,0,0,1,0)\\
	\end{array}$$
These vertices correspond to the open lattice points on the path $p(\omega_{15,\lambda})$ depicted in Figure \ref{fig:uniqueminimalnewfaceexample(15)4}.  The reader should note the position of these points relative to the lattice path $p(\omega_{15,\lambda^\star})$.  

	\begin{figure}
		\centering
		\includegraphics[width=0.7\textwidth]{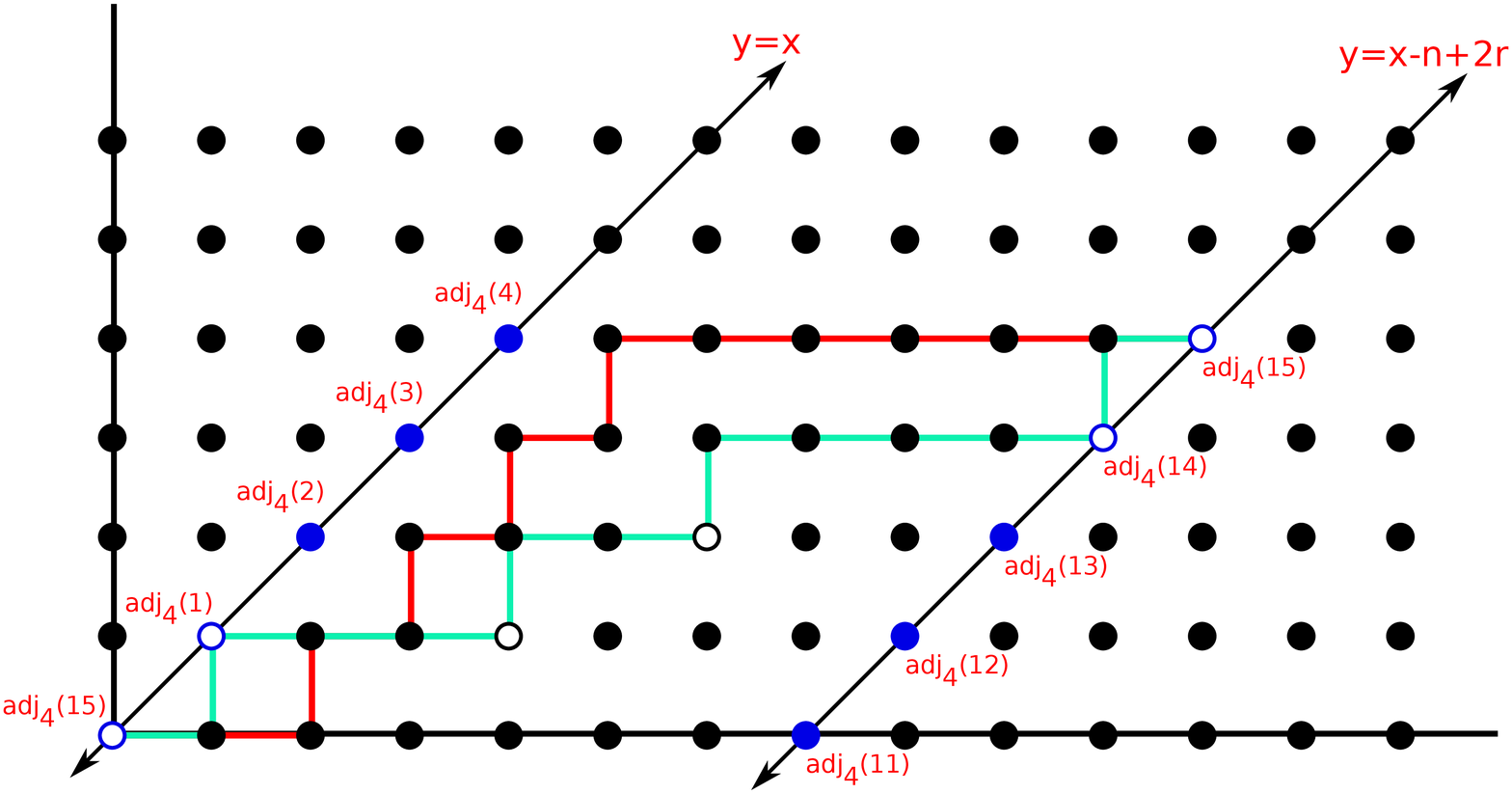}
		\caption{The lattice path $p(\omega_{15,\lambda})$ compared to $p(\omega_{15,\lambda^\star})$.  Here, $$\lambda=(1,0,0,1,0,1,0,0,0,1).$$  The unique minimal new face of $\omega_{15,\lambda}$ is given by the open lattice points.}
		\label{fig:uniqueminimalnewfaceexample(15)4}
	\end{figure}
\end{example}

To see these vertices form the unique minimal new face fix a simplex $\omega_{\ell,\lambda}\in W_{\ell,s}$, and suppose that this set of vertices is 
	$$G_{\omega}=\{v_0:=\adj{\ell},v_1,v_2,\ldots,v_q\}.$$
We will show that any face of $\omega_{\ell,\lambda}$ not using $G_\omega$ has previously appeared, and that $G_\omega$ is indeed a new face.

First consider a face $F$ of $\omega_{\ell,\lambda}$ that does not use vertex $v_t\in G_{\omega}$ for $t\neq0$.  There are then two cases:

\begin{enumerate}[(1)]
	\item $v_t\in\Adj{n}$, or
	\item $v_t\notin\Adj{n}$.
\end{enumerate}

In case (1), we saw in Corollary \ref{cortoreplacingadj} that $\omega_{\ell,\lambda}$ reduces to a previously shelled simplex that only differs for $\omega_{\ell,\lambda}$ by the vertex $v_t$.  That is, to construct a simplex $\omega_{\ell,\lambda^\prime}$ that uses the face $F$ and was shelled before $\omega_{\ell,\lambda}$ switch the parities of the moves about $v_t$.  This results in a simplex $\omega_{\ell,\lambda^\prime}\in W_{\ell,s-1}$, which was therefore shelled before $\omega_{\ell,\lambda}$.  

In case (2), we can identify a previously shelled simplex $\omega_{\ell,\lambda^\prime}\in W_{\ell,s}$ that uses the face $F$ and for which $\omega_{\ell,\lambda^\prime}<_{\colex}\omega_{\ell,\lambda}$ as follows.  Since $v_t\notin\Adj{n}$ then it must be that $v_t$ corresponds to a lattice point on $p(\omega_{\ell,\lambda})$ that is a right-most vertex on a row of the lattice, that is a corner of the path, and is also not a point on the path $p(\omega_{\ell,\lambda^\star})$ (since all points on the line $y=x$ correspond to elements of $\Adj{n}$).  Hence, the vertex $v_t$ is produced by a right move and followed by a left move in $(\omega_{\ell,\lambda})$.  Switching the parities of these moves results in replacing $v_t$ with a vertex $v_t^\prime$ which is produced by a left move and followed by a right move in the resulting cycle, say $(\omega_{\ell,\lambda^\prime})$.  Notice that the vertex $v_t^\prime$ is not an element of $\Adj{n}$.  To see this, assume otherwise.  Then by Lemma \ref{adjminnotmaxmoves} $v_t^\prime$ is a vertex produced by a left move and followed by a right move for which the number of left moves preceding the right move is maximal and not minimal with respect to equation (\ref{lambdabounds}).  Hence, the lattice point corresponding to $v_t^\prime$ lies on the line $y=x$.  But this implies that $v_t$ is a vertex on $p(\omega_{\ell,\lambda^\star})$, which is a contradiction.  Thus, $v_t^\prime$ is not an element of $\Adj{n}$.  Notice also that it is immediate from the parity switch of the moves about $v_t$ that $\omega_{\ell,\lambda^\prime}<_{\colex}\omega_{\ell,\lambda}$.  Hence, $\omega_{\ell,\lambda^\prime}\in W_{\ell,s}$ with $\omega_{\ell,\lambda^\prime}<_{\colex}\omega_{\ell,\lambda}$.  Moreover, since $\omega_{\ell,\lambda^\prime}$ uses the face $F$ since this simplex only differs from $\omega_{\ell,\lambda}$ by the vertex $v_t$, which is not used in $F$.

Now suppose $v_t=v_0=\adj{\ell}$.  By Corollary \ref{cortoreplacingadj} we know that switching the parities about $v_0$ reduces to a previously shelled simplex, which only differs from the simplex $\omega_{\ell,\lambda}$ by the vertex $v_0$.  Hence, the face $F$ also appears in the previously shelled simplex.

We next show that $G_\omega$ is indeed a new face.  Notice that by Remark \ref{latticepath} $G_\omega$ contains all the vertices in $v_{(\omega)}\cap\Adj{n}$.  For the sake of contradiction, suppose that $G_\omega$ appeared in a previously shelled simplex, say $\omega_{\ell^\prime,\lambda^\prime}$.  That is, $\omega_{\ell^\prime,\lambda^\prime}<\omega_{\ell,\lambda}$.  Since $\omega_{\ell,\lambda}\in W_{\ell,s}$ and we are assuming $\omega_{\ell^\prime,\lambda^\prime}<\omega_{\ell,\lambda}$ then $\omega_{\ell^\prime,\lambda^\prime}$ uses at most $s$ elements of $\Adj{n}$.  But since $G_\omega$ contains $s$ elements of $\Adj{n}$ we have that $\omega_{\ell^\prime,\lambda^\prime}\in W_{\ell^\prime,s}$.  In particular, $\omega_{\ell^\prime,\lambda^\prime}$ uses precisely the same elements of $\Adj{n}$ as $\omega_{\ell,\lambda}$, and so $\ell^\prime=\ell$.  Hence, $\omega_{\ell^\prime,\lambda^\prime}=\omega_{\ell,\lambda^\prime}\in W_{\ell,s}$.  So it must be that $\lambda^\prime<_{\colex}\lambda$.  That is, the right-most nonzero entry in $\lambda-\lambda^\prime$ is positive, say $\lambda_m-\lambda_m^\prime>0$.  

Consider the vertex in $\omega_{\ell,\lambda}$, say $v_t$, produced by the $m^{th}$ right move in $(\omega_{\ell,\lambda})$.  In $p(\omega_{\ell,\lambda})$ $v_t$ corresponds to the right-most corner vertex in a row of the lattice since $\lambda_m>0$.  We then have two cases.

\begin{enumerate}[(1)]
	\item The vertex $v_t$ does not correspond to a point on $p(\omega_{\ell,\lambda^\star})$.  
	\item The vertex $v_t$ does correspond to a point on $p(\omega_{\ell,\lambda^\star})$.  
\end{enumerate}

In case (1), it follows that $v_t\in G_\omega$.  Since $\omega_{\ell,\lambda}$ and $\omega_{\ell,\lambda^\prime}$ both have the same largest element in the set $\Adj{n}$, namely $\adj{\ell}$, every element of $G_\omega$ is uniquely determined by the number of left and right moves needed to produce it from $\adj{\ell}$.  In particular, we have that $v_t$ is produced by $m$ right moves and $\sum_{j=1}^{m-1}\lambda_j$ left moves in $\omega_{\ell,\lambda}$, and $v_t$ is produced by $m$ right moves and $\sum_{j=1}^{m-1}\lambda_j^\prime$ left moves in  $\omega_{\ell,\lambda^\prime}$.  But since $\lambda_m-\lambda_m^\prime>0$ and $\lambda_j-\lambda_j^\prime=0$ for all $j>m$ we have that
	$$\sum_{j=1}^{m-1}\lambda_j<\sum_{j=1}^{m-1}\lambda_j^\prime,$$
a contradiction.

\begin{figure}
	\centering
	\includegraphics[width=0.8\textwidth]{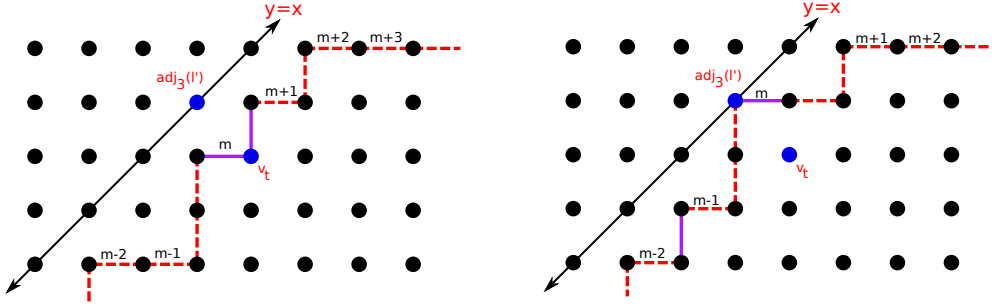}
	\caption{The case when vertex $v_t$ is a point on $p(\omega_{\ell,\lambda^\star})$.}
	\label{fig:uniqueminnewfacediagram}
\end{figure}

In case (2), it follows that $\omega_{\ell,\lambda}$ does not contain the vertex $\adj{\ell^\prime}$ corresponding to the lattice point diagonally across from $v_t$ on the line $y=x$.  This scenario is depicted if Figure \ref{fig:uniqueminnewfacediagram}.  However, since $\lambda_m-\lambda_m^\prime$ is the right-most nonzero entry in $\lambda-\lambda^\prime$ then $\omega_{\ell,\lambda^\prime}$ does contain $\adj{\ell^\prime}$, since one of the left moves accounted for by $\lambda_m$ must now be accounted for in $\lambda_t^\prime$ for $t<m$.  But this contradicts the fact that
	$$v_{(\omega_{\ell,\lambda^\prime})}\cap\Adj{n}=v_{(\omega_{\ell,\lambda})}\cap\Adj{n}.$$

Thus, $G_\omega$ is indeed a new face.  
This completes the proof of Theorem \ref{stableshellingthm} for $n$ odd.
\qed

%----THE n EVEN CASE-------------
\subsection{The $n$ even case}
\label{the n even case}
In the following, let $n$ be even.  
Following Remark~\ref{proofapproach}, our first goal is to provide a shelling of $\nabla_{n,2}^{\frac{n}{2}-1}$.  
Notice first that the labeling of the simplices in $\max\trintwor$ given by Remarks~\ref{compositionlabel} and \ref{latticepath} still apply.  
That is, $\omega\in\max\trintwor$ is uniquely described by $\omega_{\ell,\lambda}$ where $\adj{\ell}$ is the unique maximal $r$-adjacent vertex in the simplex and $\lambda = (\lambda_1,\ldots,\lambda_{n-r-1})$ is a composition of $r$ into $n-r-1$ parts for which 
$$
i-1
\leq
\sum_{j=1}^i\lambda_j
\leq
i
$$
for $i=1,\ldots,n-r-1$.  
Since $n=2r+2$ when $r=\frac{n}{2}-1$ we cannot implement the same shelling as in section~\ref{shelling} since Lemma~\ref{switchingparityforadj}, and thus Corollary~\ref{cortoreplacingadj}, no longer hold.  
In fact, switching the parity about an $r$-adjacent vertex \emph{always} produces another $r$-adjacent vertex when $n=2r+2$.  
Thus, we need to modify our shelling to accommodate for this fact.  

Recall, the order on the simplices for the $n$ odd case first collects simplices into sets $W_{\ell,s}$ where $W_{\ell,s}$ consists of all simplices with maximal $r$-adjacent vertex $\adj{\ell}$ that use $s$ elements of $\Adj{n}$.  
We then order the elements in $W_{\ell,s}$ from least-to-greatest with respect to the colexicographic order on their associated compositions.  
Finally, we order the collection of sets $\{W_{\ell,s}\}$ from least-to-greatest with respect to the colexicographic order on the labels $(\ell,s)$.  
In the following remark, we define a new order on the simplices in $\trintwor$ for $n=2r+2$ which will allow us to shell our base case.  

\begin{rmk}
\label{new order}
When $n=2r+2$ every simplex in $\max\trintwor$ uses precisely $r+1$ $r$-adjacent vertices.  
This is seen by noting that the lattice path corresponding to $\omega_{\ell,\lambda}$ lies in the region between the lines $y=x$ and $y=x-2$.  
Thus, we no longer require the parameter $s$ in the above formula.  
That is, for $n=2r+2$, we denote the collection of simplices in $\max\trintwor$ with unique maximal $r$-adjacent vertex $\adj{\ell}$ by $W_\ell$.  
We now order the elements of each set $W_\ell$ from \emph{greatest-to-least} with respect to the colexicographic order, and then order the sets $W_\ell$ as 
$$
W_1,W_2,\ldots,W_n.
$$
Denote this order on the simplices in $\max\trintwor$ by $<_e$.  
\end{rmk}

\begin{thm}
\label{gorenstein shelling theorem}
Let $n=2r+2$.  
The order $<_e$ on $\max\trintwor$ is a shelling order.
\end{thm}

\begin{proof}
To show that this order is indeed a shelling it suffices to identify the unique minimal new face of each simplex in the order.  
We claim that the unique minimal new face of a simplex $\omega$ consists of all vertices whose corresponding lattice point in $p(\omega)$ lies on the line $y=x$.  
Denote this collection by $G_\omega=\{v_0=\adj{\ell},v_1,v_2,\ldots,v_q\}$.  

We must first show that any face $F$ of $\omega$ consisting of $G_\omega\backslash\{v_t\}$ for some $t=0,\ldots,q$ appears as a face of a previously shelled simplex.  Consider first the case in which $t\neq0$.  
Let $F$ be a face of $\omega_{\ell,\lambda}$ that does not use $v_t$.  
To construct a simplex $\omega_{\ell^\prime,\lambda^\prime}<_e\omega_{\ell,\lambda}$ that uses $F$ we do the following.  
Since $t\neq0$ then $v_t$ lies on $y=x$ and switching the parities of the moves about $v_t$ replaces it with a vertex on $y=x-2$.  
This lattice path corresponds to a simplex $\omega_{\ell^\prime,\lambda^\prime}$ with $\ell^\prime=\ell$.  
Moreover, $\lambda^\prime>_{\colex}\lambda$ since switching the parities of the moves about $v_t$ amounts to shifting a move of $1$ in $\lambda$ one entry to the right.  
Thus, $\omega_{\ell,\lambda^\prime}\in W_\ell$ such that $\omega_{\ell,\lambda^\prime}<_e\omega_{\ell,\lambda}$.  

To understand the case when $t=0$ we first note that switching the parity about a point $\adj{\ell}$ replaces it with the point $\adj{\ell+r+1}$.  
This can be seen quickly from the labeling of the lattice region.  
Thus, switching the parity about $v_0$ replaces $v_0=\adj{\ell}$ with $\adj{\ell+r+1}$.  
Note $\adj{\ell+r+1}<\adj{\ell}$ whenever $\ell>r+1$.  
Since $W_1=W_2=\cdots=W_r=\emptyset$ and $W_{r+1}=\{\omega_{r+1,\lambda^\star}\}$ then $G_{\omega_{r+1,\lambda^\star}}=\emptyset$, as it is the first simplex in our shelling.  
Thus, it only remains to verify that $G_\omega$ is indeed a new face.  

To see that $G_\omega$ is a new face, we simply note that every lattice path in the decorated region corresponding to $\adj{\ell}$ for $n=2r+2$ is uniquely determined by its vertices that lie on $y=x$.  
In other words, if $G_\omega$ were indeed used previously then the simplex that used it must be $W_\ell$.  
However, the only simplex in $W_\ell$ that uses the vertices in $G_\omega$ is $\omega_{\ell,\lambda}$ itself.  
\end{proof}

\begin{cor}
\label{even stable shelling corollary}
The order $<$ defined in subsection~\ref{the shelling order} on the simplices $\omega\in\max\nabla_{n,2}^r\backslash\max\nabla_{n,2}^{r+1}$ (from least to greatest) extends the shelling of $\nabla_{n,2}^{r+1}$ to a shelling of $\nabla_{n,2}^r$.  
\end{cor}

\begin{proof}
Recall that the stable shelling of the odd second hypersimplex given in section~\ref{shelling} happens inductively, and the base case shells the smallest $(n-1)$-dimensional $r$-stable $(n,2)$-hypersimplex.  
The previous Theorem establishes the base case for $n$ even.  
Moreover, the inductive step then applies here, just as in the odd case, since the only issue was in Lemma~\ref{switchingparityforadj}.  
In the proof of this lemma, a contradiction occurs for $n=2r+2$.  
However, we are not taking this $r$ value in our inductive step.  
Thus, the same inductive shelling holds here with the Gorenstein $r$-stable second hypersimplices serving as the base case.  
\end{proof}

\begin{rmk}
\label{}
One might be concerned that to shell the Gorenstein center of the even second hypersimplices we do almost the same shelling as in the odd case, but we reverse the order of the simplices in the sets $W_\ell$.  
However, this choice is in fact consistent with our previous shelling in the sense that in both cases we always shell the simplices given by the lattice path $\lambda^\star$ first.  
\end{rmk}

Notice that Corollary~\ref{even stable shelling corollary} completes the proof of Theorem~\ref{stableshellingthm} for $n$ even.  
Thus, Corollary~\ref{even stable shelling corollary} and subsection~\ref{stable shelling in odd case proof} combine to prove Theorem~\ref{introshellingthm}.

%%%%%%%%%%%%%%%%%%%%%%%%%%%%%%%%%%%%%%%%%%%%%%%%%%%%%%%%%%%%%%%%%%%%%%%%%%%%%%%%%
%%%%%%%%%%%%%%%%%%%%%%%%%%%%%%%%%%%%%%%%%%%%%%%%%%%%%%%%%%%%%%%%%%%%%%%%%%%%%%%%%

\subsection{Some Corollaries to Theorem~\ref{introshellingthm}}

We first note that inductively Theorem \ref{stableshellingthm} results in a shelling of the second hypersimplex $\hyperntwo$, thereby proving Theorem \ref{introshellingthm}.  
This shelling is interesting in the sense that it begins with simplices that use only the ``most stable" vertices of the polytope and at each stage adds a simplex that uses more and more of the ``less stable" vertices.  
We now give a few results that will be helpful in section \ref{hstarpolynomials}, where we examine the $h^\ast$-polynomials of these polytopes.

\begin{cor}\label{maxdimminnewface}
Let $\omega_{\ell,\lambda}\in\max\nabla_{n,2}^r\backslash\max\nabla_{n,2}^{r+1}$.  The maximum dimension of the minimum new face $G_\omega$ is $r+1$.
\end{cor}

\begin{proof}
Consider the lattice path $p(\omega_{\ell,\lambda})$ and recall that the vertices of the minimal new face $G_\omega$ correspond to lattice points on $p(\omega_{\ell,\lambda})$ that lie on the line $y=x$ together with those that are the right-most such points in their row of the lattice, that are corners of $p(\omega_{\ell,\lambda})$, and are not on the path $p(\omega_{\ell,\lambda^\star})$.  In particular, these are the lattice points $(x,x)$ for $x\in[r]$ (which we will call \emph{type 1}), and $(x,y)$  where $y\leq x-3$, $y\in\{0,1,\ldots,r-1\}$, and $(x,y)$ is a corner of $p(\omega_{\ell,\lambda})$ that is not also on the path $p(\omega_{\ell,\lambda^\star})$ (which we will call \emph{type 2}).  The lattice path $p(\omega_{\ell,\lambda})$ can have at most one such point  on each line $y=\alpha$ for $\alpha\in[r]$, and at most two such points on the line $y=0$ (one of which is always $\adj{\ell}$).  This gives a maximum possible dimension of $r+2$ for $G_\omega$.  

Assume now that there exists $\omega_{\ell,\lambda}$ such that $G_{\omega}$ has dimension $r+2$.    Consider the vertices of $G_\omega$ corresponding to lattice points $(x,x)$ for $x\in[r]$.  Label this set of vertices by $V$, and for $v\in V$ label its corresponding lattice point by $p_v:=(x_v,y_v)$.  Notice that $V$ is nonempty since by our assumption $G_\omega$ contains a vertex corresponding to a lattice point on the line $y=r$, and (since we do not wish to over-count the vertex $\adj{\ell}$) the only option is $\adj{\ell+r}\in V$.  So let $v\in V$.  Then $p(\omega_{\ell,\lambda})$ cannot contain any type 2 points on the lines $y=y_v-1$ or $y=y_v-2$ (this is because $p(\omega_{\ell,\lambda})$ may only use North and East moves).  But by our assumption the lines $y=y_v-1$ and $y=y_v-2$ must each contain a point corresponding to an element of $G_\omega$.  Hence, they are type 1 points, and therefore elements of $ V$.  

Beginning with $v=\adj{\ell+r}$, which has corresponding lattice point $p_v=(r,r)$, iterating this argument shows that for each line $y=\alpha$, $\alpha\in\{0,1,2,\ldots,r-1\}$, the path $p(\omega_{\ell,\lambda})$ uses the point $(\alpha,\alpha)$, and no type 2 points on $y=\alpha$.  Hence, $\#(G_\omega)=r+1$, a contradiction.
\end{proof}

For the case when $r=1$, the next corollary is also a corollary to the algebraic formula given for the $h^\ast$-polynomial of $\Delta_{n,2}$ by Katzman in \cite{kat}.  However, we are now able to give an entirely combinatorial proof of this result.

\begin{cor}\label{degreeofhstar}
For $r<\left\lfloor\frac{n}{2}\right\rfloor$ the degree of the $h^\ast$-polynomial of $\rstablentwo$  is $\left\lfloor\frac{n}{2}\right\rfloor$.  
When $n$ is odd the leading coefficient is $n$, and when $n$ is even the leading coefficient is $1$.
\end{cor}

\begin{proof}
When $n$ is even the result follows immediately from Theorem~\ref{gorenstein shelling theorem} and Corollary~\ref{maxdimminnewface}.  
Since $r<\left\lfloor\frac{n}{2}\right\rfloor$ we have shelled the simplices in $\max\nabla_{n,2}^{\left\lfloor\frac{n}{2}\right\rfloor-1}\backslash\max\nabla_{n,2}^{\left\lfloor\frac{n}{2}\right\rfloor}$ in order to build the hypersimplex $\rstablentwo$.  By Corollary \ref{maxdimminnewface} for a simplex $\omega\in\max\nabla_{n,2}^{\left\lfloor\frac{n}{2}\right\rfloor-1}\backslash\max\nabla_{n,2}^{\left\lfloor\frac{n}{2}\right\rfloor}$ the maximum dimension of $G_\omega$ is $\left\lfloor\frac{n}{2}\right\rfloor$.  It remains to show that this maximum dimension is achieved precisely $n$ times.  

Notice first that for $r=\left\lfloor\frac{n}{2}\right\rfloor-1$ we have that the lattice paths labeling simplices in \\
$\max\nabla_{n,2}^r\backslash\max\nabla_{n,2}^{r+1}$ are bounded between the lines $y=x$ and $y=x-3$.  
Also, for $\omega$ to satisfy $\#(G_{\omega})=r+1$ then there must be a total of $r$ points of $p(\omega_{\ell,\lambda})$ on the lines $y=x$ and $y=x-3$ other than $(0,0)$ and $(n-r,r)$.  

For $\max(v_{(\omega)}\cap\Adj{n})=\adj{\ell}$ with $\ell\leq r$ this is impossible since there are less than $r$ points on these lines that we may use without violating the choice of $\max(v_{(\omega)}\cap\Adj{n})=\adj{\ell}$.  

For $r<\ell<n$ consider the following.  Suppose $p(\omega_{\ell,\lambda})$ uses a point $(\alpha,\alpha)$ for $0<\alpha\leq r$.  Then by the same argument as in Corollary \ref{maxdimminnewface} this implies that $p(\omega_{\ell,\lambda})$ uses $(\alpha-1,\alpha-1)$.  Iterating this just as before we get that $p(\omega_{\ell,\lambda})$ uses the points 
	$$\{(0,0),(1,1),(2,2),\ldots,(\alpha,\alpha)\}.$$
However, this contradicts the fact that $r\leq\ell<n$ and $\max(v_{(\omega)}\cap\Adj{n})=\adj{\ell}$.  Hence, the only point used by $p(\omega_{\ell,\lambda})$ on the line $y=x$ is $(0,0)$.  Therefore, $p(\omega_{\ell,\lambda})$ must be the path using all the points on the line $y=x-3$.

For $\ell=n$ there are exactly $r+1$ paths such that $\#(G_\omega)=r+1$.  This is also seen from the iterative argument we used in Corollary \ref{maxdimminnewface}.  Suppose the path $p(\omega_{\ell,\lambda})$ uses $s\leq r$ points of the set $\{(1,1),(2,2),\ldots,(r,r)\}$, and let $(\alpha,\alpha)$ be the point in this collection for which the value of $\alpha$ is maximal.  It then follows that $p(\omega_{\ell,\lambda})$ uses all the points 
	$$\{(0,0),(1,1),(2,2),\ldots,(\alpha,\alpha)\}.$$
Hence, it must be that $\alpha=s$.  Since there is exactly one path that uses $r+1$ points on the lines $y=x$ and $y=x-3$ and uses the points $\{(0,0),(1,1),(2,2),\ldots,(\alpha,\alpha)\}$, then we conclude that there are exactly $r+1$ simplices $\omega_{n,\lambda}$ with $\#(G_\omega)=r+1$.  

Considering all of these cases together we conclude that there are exactly $n$ simplices with $\#(G_\omega)=r+1$.  
\end{proof}

\begin{rmk}
We remark that the proofs of the previous corollaries are intriguing since they point out that this shelling allows us to study the Ehrhart theory of the second hypersimplices, as well as the $r$-stable second hypersimplices,  by enumerating lattice paths in various ladder-shaped regions of the plane.  However, this enumeration problem, in general, is not trivial as suggested by the work of Krattenthaler in \cite{krat}.
\end{rmk}

%%%%%%%%%%%%%%%%%%%%%%%%%%%%%%%%%%%%%%%%%%%%%%%%%%%%%%%%%%%%%%%%%%%%%%%%%%%%%%%%%%
%%%%%%%%%%%%%%%%%%%%%%%%%%%%%%%%%%%%%%%%%%%%%%%%%%%%%%%%%%%%%%%%%%%%%%%%%%%%%%%%%%

\section{$h^\ast$-polynomials}\label{hstarpolynomials}

In the following we compute the $h^*$-polynomial of $\rstablentwo$.
For every $1\leq r\leq \floorntwo$ we give a formula for the $h^*$-polynomial of $\rstablentwo$ in terms of a sum of independence polynomials of certain graphs.
In the case of $n$ even and $r=\floorntwo-1$ we show that the $h^*$-polynomial of $\rstablentwo$ is the generating polynomial for the binomial coefficients ${r+1\choose m}$.  
In the case of $n$ odd and $r=\floorntwo-1$ we show that the $h^*$-polynomial of $\rstablentwo$ is precisely the independence polynomial of the cycle on $n$ vertices, or equivalently, the $n^{th}$ Lucas polynomial.  
Then, via Ehrhart-MacDonald reciprocity, we demonstrate that the $h^*$-polynomial of the relative interior of this polytope is a univariate specialization of a polynomial that plays an important role in the theory of proper holomorphic mappings of complex balls in Euclidean space.
Specifically, this polynomial is the squared Euclidean norm function of a well-studied CR mapping of the Lens space into the unit sphere within $\C^{r+3}$ \cite{dangelo1, dangelo2, dangelo3}.
We end with a discussion of when these $h^\ast$-polynomials appear to be unimodal and make two conjectures along these lines.  

%More generally, for $r\leq \floorntwo-1$, we show that the $h^\ast$-polynomial of the relative interior of $\rstablentwo$ is a univariate specialization of the squared Euclidean norm function of certain polynomial maps that induce smooth immersions of the Lens Space $L(n,2)$ into $\C^{n+1}$.  
%In the case that $r=\floorntwo-1$, this immersion is a well-studied CR mapping of the Lens space into the sphere \cite{dangelo1, dangelo2, dangelo3}.

%%%%%%%%%%%%%%%%%%%%%%%%%%%%%%%%%%%%%%%%%%%%%%%%%%%%%%%%%%
%%%%%%%%%%%%%%%%%%%%%%%%%%%%%%%%%%%%%%%%%%%%%%%%%%%%%%%%%%

\subsection{The $h^*$-polynomial of $\rstablentwo$ via Independence Polynomials of Graphs.}\label{independencepolynomials}

It will be helpful to recall some basic facts about independence polynomials of graphs.
Suppose that $G$ is a finite simple graph with vertex set $V(G)$ and edge set $E(G)$.  
An \emph{independent} set in $G$ is a subset of the vertices of $G$, $S\subset V(G)$, such that no two vertices in $S$ are adjacent in $G$.  
Let $s_i$ denote the number of independent sets in $G$ with cardinality $i$, and let $\beta(G)$ denote the maximal size of an independent set in $G$.  
The \emph{independence polynomial} of $G$ is the polynomial
	$$I(G;x):=\sum_{i=0}^{\beta(G)}s_ix^i.$$
Independence polynomials of graphs are well-studied structures.  
Levit and Mandrescu nicely survey properties of these polynomials in \cite{levit1}.  
Here, we restrict our attention only to those properties which we will use to compute $h^*$-polynomials.
Suppose that $G_1$ and $G_2$ are two finite simple graphs, and let $G_1\cup G_2$ denote their disjoint union.  
It is a well known fact (see for instance \cite{gutman} or \cite{levit1}) that 
	$$I(G_1\cup G_2;x)=I(G_1;x)\cdot I(G_2;x).$$
Let $P_n$ and $C_n$ denote the path and cycle on $n$ vertices, respectively.  
In \cite{arocha}, Arocha showed that
	$$I(P_n;x)=F_{n+1}(x) \qquad \mbox{and} \qquad I(C_n;x)=F_{n-1}(x)+2xF_{n-2}(x),$$
where $F_n(x)$ denotes the $n^{th}$ \emph{Fibonacci polynomial}.  
The Fibonacci polynomials are defined for $n\geq0$ by the recursion
	$$F_0(x)=1,\quad F_1(x)=1,\quad \mbox{and} \quad F_n(x)=F_{n-1}(x)+xF_{n-2}(x).$$
A closely related class of polynomials are the \emph{Lucas polynomials}, which are defined by the recursion
	$$L_0(x)=2,\quad L_1(x)=1,\quad \mbox{and} \quad L_n(x)=L_{n-1}(x)+xL_{n-2}(x).$$
These collections of polynomials will play important roles in our computations of $h^*$-polynomials.

In the following we let $\hstar{\rstablentwo}$ denote the $h^*$-polynomial of $\rstablentwo$.  
Recall that Theorem \ref{stableshellingthm} provides a shelling of the unimodular triangulation of $\rstablentwo$ induced by the circuit triangulation of $\hyperntwo$.  
We let $\trintwor$ denote this triangulation of $\rstablentwo$ and $\max\trintwor$ denote the collection of maximal simplices in $\trintwor$.
By a theorem of Stanley \cite{stan3}, we may compute 
	$$\hstar{\rstablentwo}=\sum_{i=0}^{n-1}h_i^*x^i$$
where $h_i^*$ equals the number of simplices in $\max\trintwor$ with unique minimal new face of dimension $i-1$ with respect to the shelling described by Theorem \ref{stableshellingthm}.  
Also recall that by Lemma \ref{fulldim}
	$$\left\{\Delta_{n,2}^{stab\left(\left\lfloor\frac{n}{2}\right\rfloor\right)}\right\}=\max\nabla_{n,2}^{\left\lfloor\frac{n}{2}\right\rfloor}$$
for $n$ odd, and this simplex has unique minimal new face $\emptyset$.  
Moreover, every simplex in $\max\trintwor$ other than this simplex has unique minimal new face of dimension at least 0.  
Analogously, for $n$ even, every simplex in $\max\trintwor$ has unique minimal new face of dimension at least $0$ except for $\omega_{r+1,\lambda^\star}$ for $r=\frac{n}{2}-1$. 
For a fixed $1\leq r\leq\floorntwo$, the simplices in $\max\trintwor\backslash\max\nabla_{n,2}^{r+1}$ correspond to lattice paths in decorated ladder-shaped regions of the plane.
Since the shape of this region is fixed for fixed values of $n$ and $r$, and there is one such region for each $\ell\in[n]$, we will refer to the decorated ladder-shaped region with origin label $\ell$ as an \emph{$\ell$-region}.

Suppose $\omega\in\max\trintwor\backslash\max\trintworplus$ is a simplex with corresponding lattice path $\lambda(\omega)$.  
Then, for a fixed $\ell$, $\lambda(\omega)$ is a lattice path in the corresponding $\ell$-region if and only if $\adj{\ell}$ is the maximal $r$-adjacent vertex used by $\omega$.
Hence, any point on the boundary lines $y=x$ and $y=x-n+2r$ labeled by $s$ with $s>\ell$ cannot be used by $\lambda(\omega)$.
In this way, these lattice points are \emph{inaccessible} lattice points of the $\ell$-region.  
We now formalize these definitions.

\begin{defn}
Fix $n>2$ and $1\leq r\leq\floorntwo$.
For each $\ell\in[n]$, we call the decorated ladder-shaped region of $\ZZ^2$ containing the lattice paths corresponding to simplices in $\max\trintwor$ with maximal $r$-adjacent vertex $\adj{\ell}$, an {\bf $\ell$-region}.
For a fixed $\ell\in[n]$, an {\bf accessible} point in an $\ell$-region is a lattice point that does not lie on the path $\lambda^\star$, but does lie on a path $\lambda(\omega)$ for some simplex $\omega\in\max\trintwor$ with maximal $r$-adjacent vertex $\adj{\ell}$.
%For a fixed $\ell\in[n]$, we call a lattice point of an $\ell$-region that lies on $y=x$ or $y=x-n+2r$ {\bf accessible} if it has label $s$ with $s<\ell$.  
Otherwise, it is called {\bf inaccessible}.
\end{defn}

Recall that the vertices of the unique minimal new face of $\omega$ correspond to the lattice points labeled by $\ell$ and the corners of $\lambda(\omega)$ that ``point away" from the path $\lambda^\star$.
With these facts in hand, we are ready to compute the $h^*$-polynomials of $\rstablentwo$ for $n$ odd.

We first examine the case where $n$ is odd and $r=\floorntwo-1$.  
Then we will extend this result to the remaining values of $n$ and $r$.
For $n$ is odd and $r=\floorntwo-1$, the boundary lines of the $\ell$-regions are given by $y=x$ and $y=x-3$.  
Since the vertices of the unique minimal new face of $\omega$ correspond to the lattice points labeled by $\ell$ and the corners of $\lambda(\omega)$ that ``point away" from the path $\lambda^\star$ then the dimension of the unique minimal new face of $\omega$ is the number of corners of $\lambda(\omega)$ that lie on the lines $y=x$ and $y=x-3$.  
(Notice that the \emph{cardinality} of the unique minimal new face is one more than this value, in which case we include the vertex corresponding to the origin in our count as well.)
Suppose that $(x,x)$ is an accessible lattice point in an $\ell$-region.  
Then the lattice path $\lambda(\omega)$ may either use the lattice point in the set $\{(x,x)\}$ or it may use  lattice points in the set $\{(x+2,x-1),(x+1,x-2)\}$, but it may not use points from both sets.
This is an immediate consequence of the fact that $\lambda(\omega)$ only uses North and East moves.
In this way, the lattice point $(x,x)$ \emph{inhibits} the lattice points in the set $\{(x+2,x-1),(x+1,x-2)\}$ and vice versa.
For convenience, we make this a formal definition.

\begin{defn}
Fix $n>2$, $1\leq r\leq \floorntwo$, and $\ell\in[n]$.  
Suppose $a$ and $b$ are lattice points in the corresponding $\ell$-region.  
We say that $a$ and $b$ {\bf inhibit} one another if and only if the vertices corresponding to $a$ and $b$ cannot appear together in the unique minimal new face of any simplex in $\max\trintwor\backslash\max\trintworplus$ with maximum $r$-adjacent vertex $\adj{\ell}$.  
Moreover, if $a$ and $b$ possess labels $\ell_1$ and $\ell_2$, we also say that $\ell_1$ and $\ell_2$ {\bf inhibit} one another.
\end{defn}

Using the labels of these points in the $\ell$-region this scenario can be represented in the following manner.  
Consider a line of $2r$ spots labeled as follows:

\begin{center}
	\includegraphics[width=0.7\textwidth]{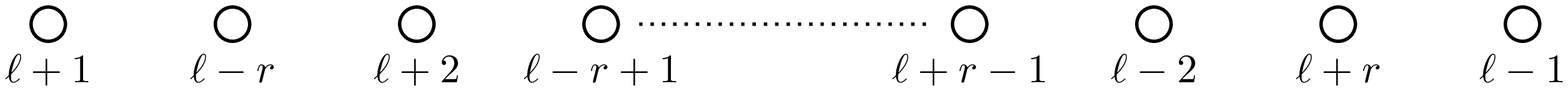}
\end{center}

Notice that each label is adjacent to exactly those labels that it inhibits.  
For each $\ell\in[n]$ we will use this diagram to construct a graph $G_{n,r,\ell}$ whose independent sets (together with the origin) are precisely the unique minimal new faces of the simplices whose lattice paths reside in the $\ell$-region.  
Fix $\ell\in[n]$, and let $S$ denote the set of accessible vertices in the corresponding $\ell$-region. 
Construct the graph $G_{n,r,\ell}$ by filling in each spot in the above diagram corresponding to a vertex in $S$.  We think of these filled in spots as the vertices of $G_{n,r,\ell}$, and we place an edge between any two vertices that are not separated by a spot.

\begin{example}\label{lineinhibitiondiagramexample}
Here are the graphs $G_{n,r,\ell}$ for each choice of $\ell$ when $n=13$ and $r=\floorntwo-1=5$.

	\begin{center}
		\includegraphics[width=0.7\textwidth]{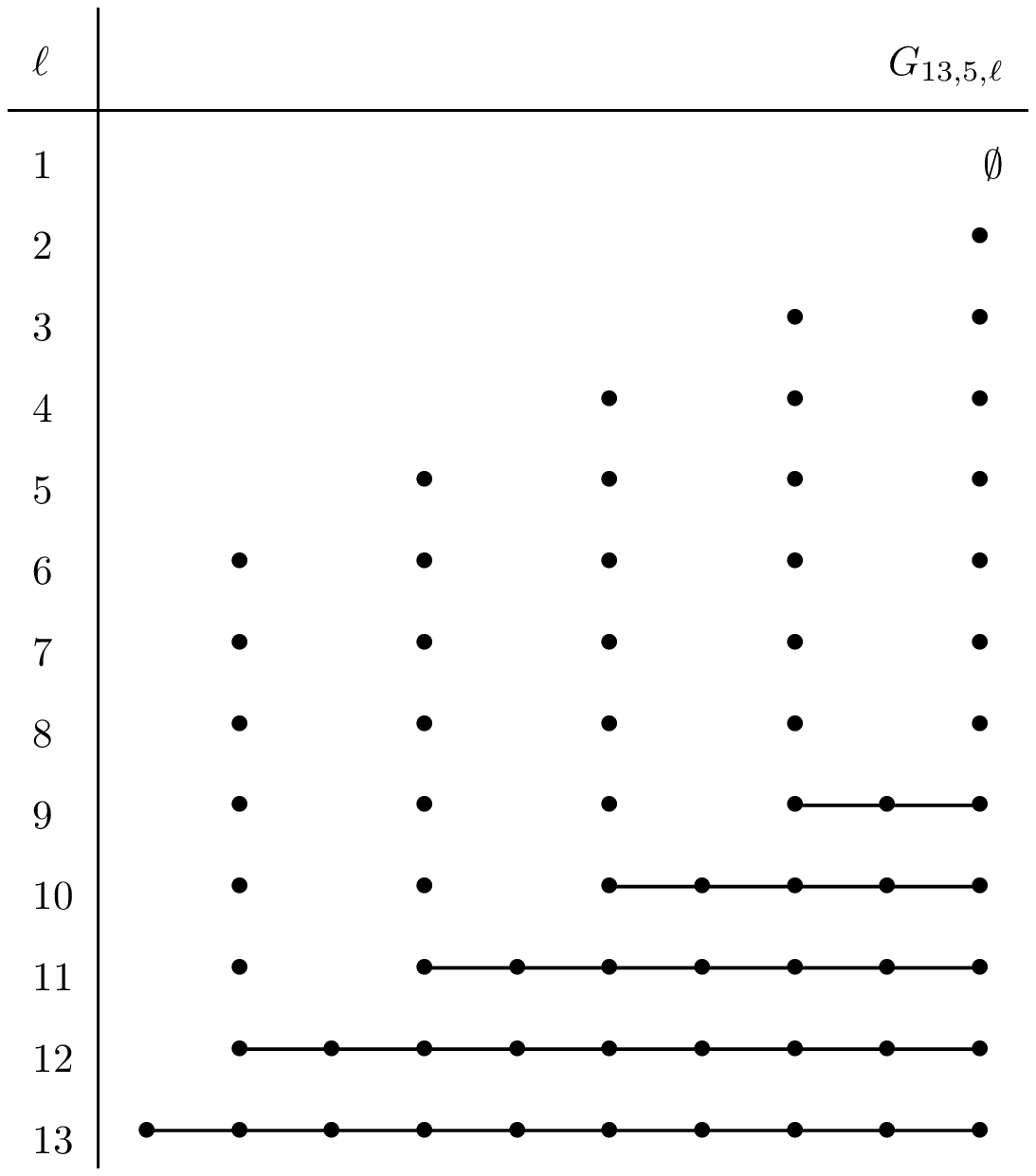}
	\end{center}
\end{example}

\begin{prop}\label{sumformula}
Fix odd $n>2$ and let $r=\floorntwo-1$.  Then 
	$$\hstar{\rstablentwo}=1+x\sum_{\ell=1}^nI(G_{n,r,\ell};x).$$
Equivalently,
	$$\hstar{\rstablentwo}=1+x\left(\sum_{i=0}^{r-1}(1+x)^i+3(1+x)^r+\sum_{i=0}^{r-2}(1+x)^iF_{2r-2i}(x)+F_{2r+1}(x)\right).$$
\end{prop}

\begin{proof}
Suppose 
	$$I(G_{n,r,\ell};x)=s_0+s_1x+s_2x^2+\cdots+s_{\beta(G_{n,r,\ell})}x^{\beta(G_{n,r,\ell})}.$$
Then $s_i$ is equal to the number of independent sets of cardinality $i$ in $G_{n,r,\ell}$.
By the construction of $G_{n,r,\ell}$ this number is precisely the number of simplices in $\max\trintwor\backslash\max\trintworplus$ that have unique minimal new face of dimension $i$ and maximal $r$-adjacent vertex $\adj{\ell}$.  Hence,
$$\hstar{\rstablentwo}=1+x\sum_{\ell=1}^nI(G_{n,r,\ell};x).$$

To prove the second equality we must identify the graphs $G_{n,r,\ell}$.  
In order to do this, we must understand the accessible lattice points in an $\ell$-region for each $\ell\in[n]$.  
In general, the set of accessible points for each $\ell$-region is given by

\begin{center}
\begin{tabular}{c|l}
$\ell$ & Set of accessible points in the $\ell$-region \\
\hline
$1$ & $\{\ell\}$ \\
$2$ & $\{\ell,\ell-1\}$ \\
$3$ & $\{\ell,\ell-1,\ell-2\}$ \\
\vdots & \vdots \\
$r$ & $\{\ell,\ell-1,\ell-2,\cdots,\ell-r+1\}$ \\
$r+1$ & $\{\ell,\ell-1,\ell-2,\cdots,\ell-r+1,\ell-r\}$ \\
$r+2$ & $\{\ell,\ell-1,\ell-2,\cdots,\ell-r+1,\ell-r\}$ \\
$r+3$ & $\{\ell,\ell-1,\ell-2,\cdots,\ell-r+1,\ell-r\}$ \\
$r+4$ & $\{\ell,\ell-1,\ell-2,\cdots,\ell-r+1,\ell-r,\ell+r\}$ \\
$r+5$ & $\{\ell,\ell-1,\ell-2,\cdots,\ell-r+1,\ell-r,\ell+r,\ell+r-1\}$ \\
$r+6$ & $\{\ell,\ell-1,\ell-2,\cdots,\ell-r+1,\ell-r,\ell+r,\ell+r-1,\ell+r-2\}$ \\
\vdots & \vdots \\
$n-1$ & $\{\ell,\ell-1,\ell-2,\cdots,\ell-r+1,\ell-r,\ell+r,\ell+r-1,\ell+r-2,\ldots,\ell+2\}$ \\
$n-1$ & $\{\ell,\ell-1,\ell-2,\cdots,\ell-r+1,\ell-r,\ell+r,\ell+r-1,\ell+r-2,\ldots,\ell+2,\ell+1\}$ \\
\end{tabular}
\end{center}

That is, as $\ell$ increases, we first gain points on the diagonal $y=x-3$ from top-to-bottom, and then those on $y=x$ from top-to-bottom.  
This happens one point at a time except for $\ell=r+2$ and $\ell=r+3$, which have the same set of accessible points as $\ell=r+1$.  
This is because $n=2r+3$ when $r=\floorntwo-1$, and the labels on the diagonals of an $\ell$-region are precisely those vertices of a convex $n$-gon labeled $1,2,\ldots,n$ with circular distance at most $r$ from $\ell$.

Given this characterization of the accessible points in each $\ell$-region we can determine the graphs $G_{n,r,\ell}$ as follows.  
For $\ell\leq r$ the graph $G_{n,r,\ell}$ is simply a collection of $\ell-1$ disjoint vertices.  
For $\ell\in\{r+1,r+2,r+3\}$ the graph $G_{n,r,\ell}$ is a collection of $r$ disjoint vertices.  
For $r+4\leq \ell\leq n-1$ we begin to add edges.  
For convenience let $\ell=r+3+t$ for the suitable value of $t$.  
Then the graph $G_{n,r,\ell}$ includes the vertices 
	$$\ell+r,\ell+r-1,\ell+r-2,\ldots,\ell+r-(t-1).$$
Each of these vertices attaches to each vertex next to it in the diagram.  
This produces a path of length $2t+1$ and a collection of $r-(t+1)$ disjoint vertices.  
Finally, for $\ell=n$, all points are accessible so $G_{n,r,\ell}=P_{2r}$.
Then, if we make the change of variables $i=r-(t+1)$, we have that

\begin{equation*}
\begin{split}
\hstar{\rstablentwo}&=1+x\sum_{\ell=1}^nI(G_{n,r,\ell};x),\\
&=1+x\left(\sum_{i=0}^{r-1}(1+x)^i+3(1+x)^r+\sum_{i=0}^{r-2}(1+x)^iI(P_{2r-2i-1};x)+I(P_{2r};x)\right),\\
&=1+x\left(\sum_{i=0}^{r-1}(1+x)^i+3(1+x)^r+\sum_{i=0}^{r-2}(1+x)^iF_{2r-2i}(x)+F_{2r+1}(x)\right).\\
\end{split}
\end{equation*}
\end{proof}

This formula for $\hstar{\rstablentwo}$ is convenient since it allows us to compute this polynomial via rows and diagonals of Pascal's Triangle.
In subsection \ref{lens}, we will see that this formula is equal to the $n^{th}$ Lucas polynomial.
However, we first show how we may generalize this formula to the remaining $r$-stable $(n,2)$-hypersimplices.  
We begin with a proposition analogous to Proposition~\ref{sumformula} in the case when $n=2r+2$, thereby providing a formula for $\hstar{\rstablentwo}$ for $r=\floorntwo-1$ for every $n$.  
\begin{prop}
\label{even gorenstein hstarpolynomials theorem}
Let $n=2r+2$.  
Then 
$$
\hstar{\rstablentwo} = (x+1)^{r+1}.
$$
\end{prop}

\begin{proof}
Using the notation from subsection~\ref{the n even case}, we count the simplices contributing to the coefficient $h^\ast_i$ by counting those arising from each set $W_\ell$, $\ell=1,\ldots,n$.  
For $\ell = 1,\ldots,r$, we have that $W_\ell=\emptyset$, so nothing is contributed the coefficient $h^\ast_i$.  
For $\ell=r+1$, we have that $W_{\ell}=\{\omega_{\ell,\lambda^\star}\}$, and this simplex has the unique minimal new face $\emptyset$.  
Thus, $h_0=1$.  

For $\ell>r+1$, we have that $\ell=r+2+t$ for some $0\leq t\leq r$.  
Moreover, the lattice point on the line $y=x$ in the ladder shaped region assoicated to $\adj{\ell}$ that is furthest from the origin is $\adj{\ell+r}=\adj{2r+2+t}=\adj{t}$.  
Thus, for $\ell=r+2+t$, there are $t$ available lattice points on $y=x$ (other than the origin).  
Since any choice of $s\leq t$ of these available points corresponds to a unique path in the region, which in turn corresponds to a unique simplex in $\max\trintwor$ with unique minimal new face of dimension $s-1$, then the polynomial
$$
h_1^\ast x+h_2^\ast x^2+\cdots+h_{n-1}^\ast x^{n-1} = x(x+1)^t.
$$
Thus, 
$$
\hstar{\rstablentwo} =
1+x\sum_{t=0}^r(x+1)^t =
(x+1)^{r+1}.
$$
\end{proof}

We have now established a formula for $\hstar{\rstablentwo}$ for $r=\floorntwo-1$ for all $n$.  
To extend these formulas to $\hstar{\rstablentwo}$ for $1\leq r<\floorntwo-1$, we first construct an inhibition diagram similar to the one used for $r=\floorntwo-1$.  

Fix $1\leq r<\floorntwo-1$.
Construct a graph $G_{n,r}$ with vertex set consisting of all lattice points in the $n$-region that do not lie on the path $\lambda^\star$, and edge set $\{\{i,j\}: \mbox{$i$ inhibits $j$}\}$.  
Let $G_{n,r,\ell}$ denote the subgraph of $G_{n,r}$ that is induced by the set of accessible points in the $\ell$-region.
We refer to the graph $G_{n,r,\ell}$ as the \emph{inhibition diagram} for its associated $\ell$-region.

\begin{prop}\label{generalsumformula}
Fix $n>2$ and $1\leq r\leq \floorntwo-1$.  
Then
	$$\hstar{\rstablentwo}=1+x\left(\sum_{j=r}^{\floorntwo-1}\sum_{\ell=1}^nI(G_{n,j,\ell};x)\right).$$
In particular, the $h^\ast$-polynomial of $\hyperntwo$ is given by a sum of independence polynomials.
\end{prop}

\begin{proof}
The proof of this formula is identical to the proof of the formulas given in Propositions~\ref{sumformula} and \ref{even gorenstein hstarpolynomials theorem}.
\end{proof}

Notice that Proposition \ref{generalsumformula} proves Theorem \ref{introsumsofindependencepolysthm}.  
We end this subsection with a few examples of the graphs $G_{n,r,\ell}$ and the resulting $h^*$-polynomials.

\begin{example}
\label{inhibitiondiagrams}
In Example \ref{lineinhibitiondiagramexample} we saw that the inhibition diagrams for $\ell$-regions when $r=\floorntwo-1$ correspond to subgraphs of a path of length $2r$.  
Suppose now that $r=\floorntwo-2$.  
Here are the general $\ell$-regions for $n=7,9$, and $11$, respectively.
\begin{center}
$\begin{array}{ccc}
		\includegraphics[width=0.3\textwidth]{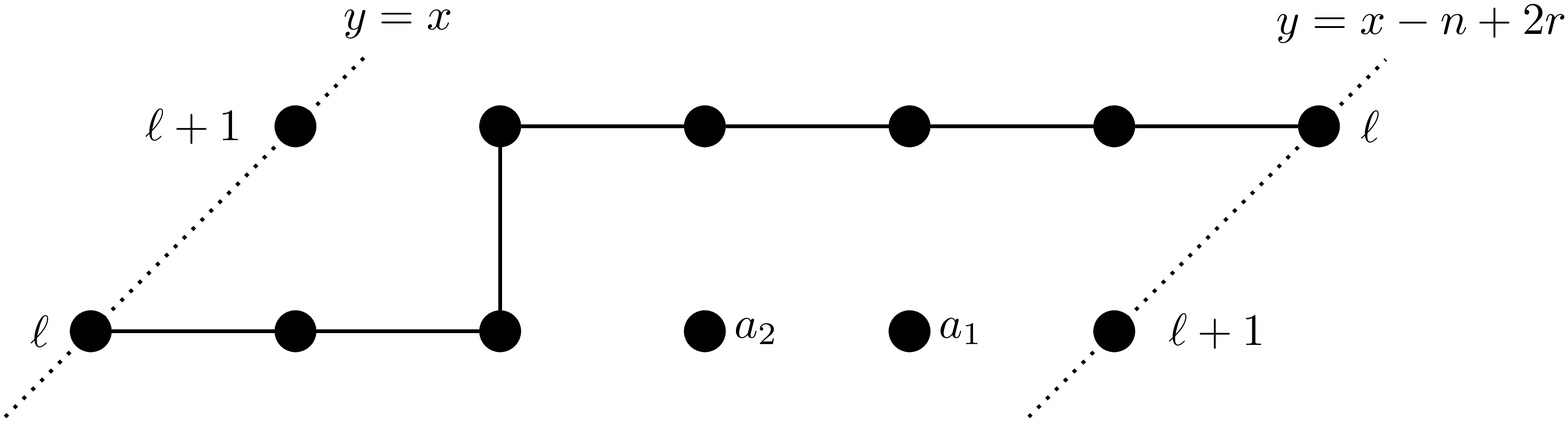} & \includegraphics[width=0.3\textwidth]{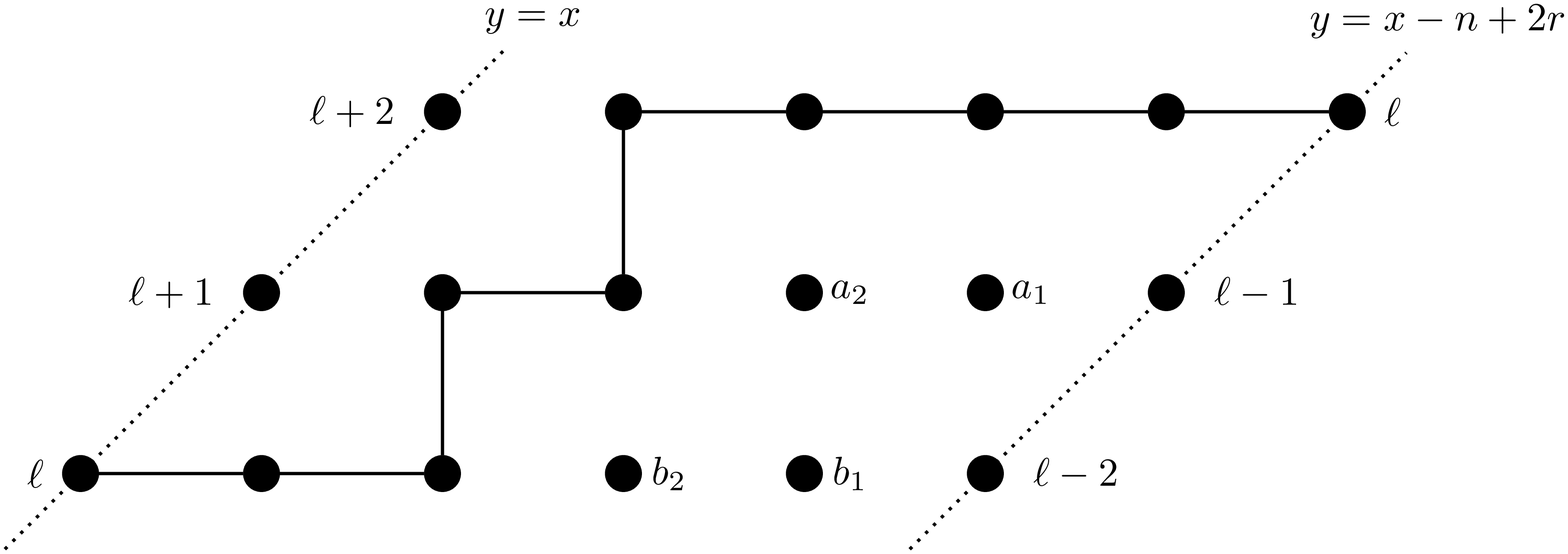} & \includegraphics[width=0.3\textwidth]{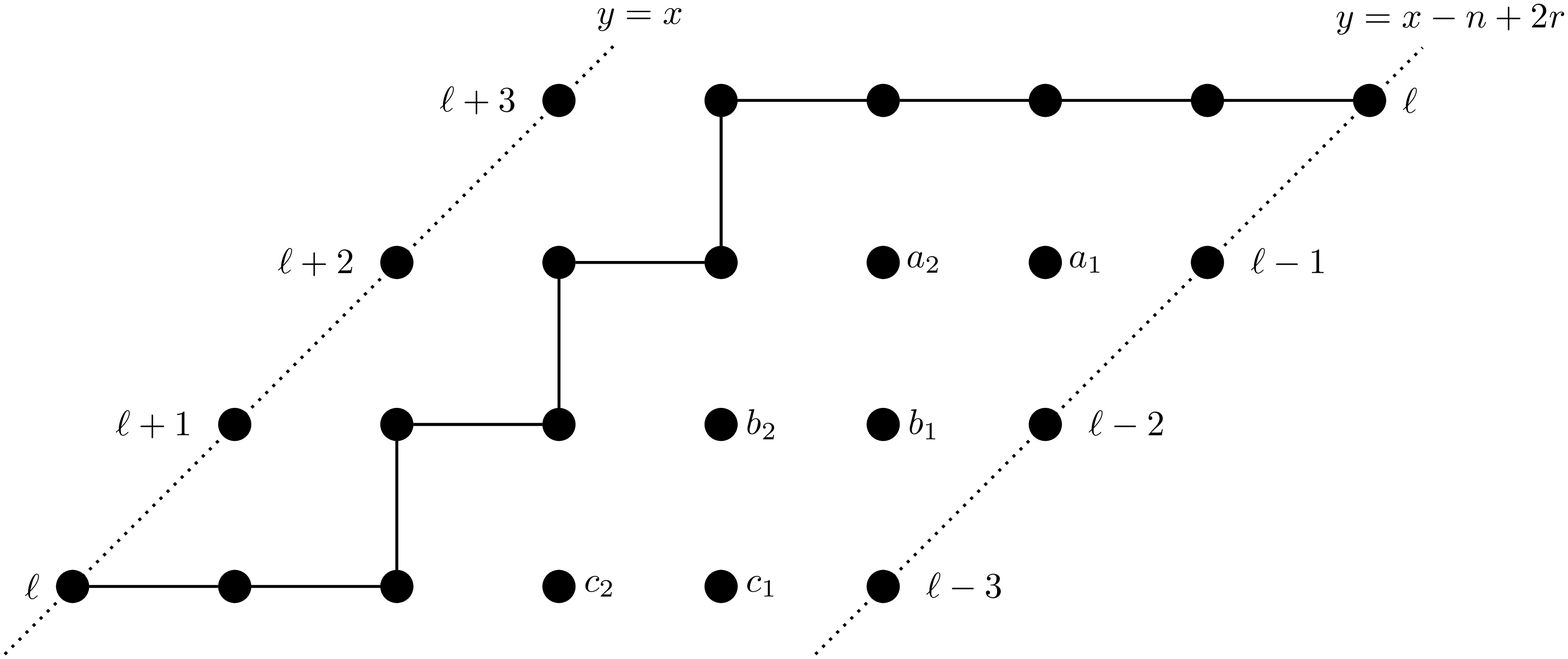}
\end{array}$
\end{center}
The corresponding inhibition diagrams are, respectively,
\begin{center}
$\begin{array}{ccc}
		\includegraphics[width=0.3\textwidth]{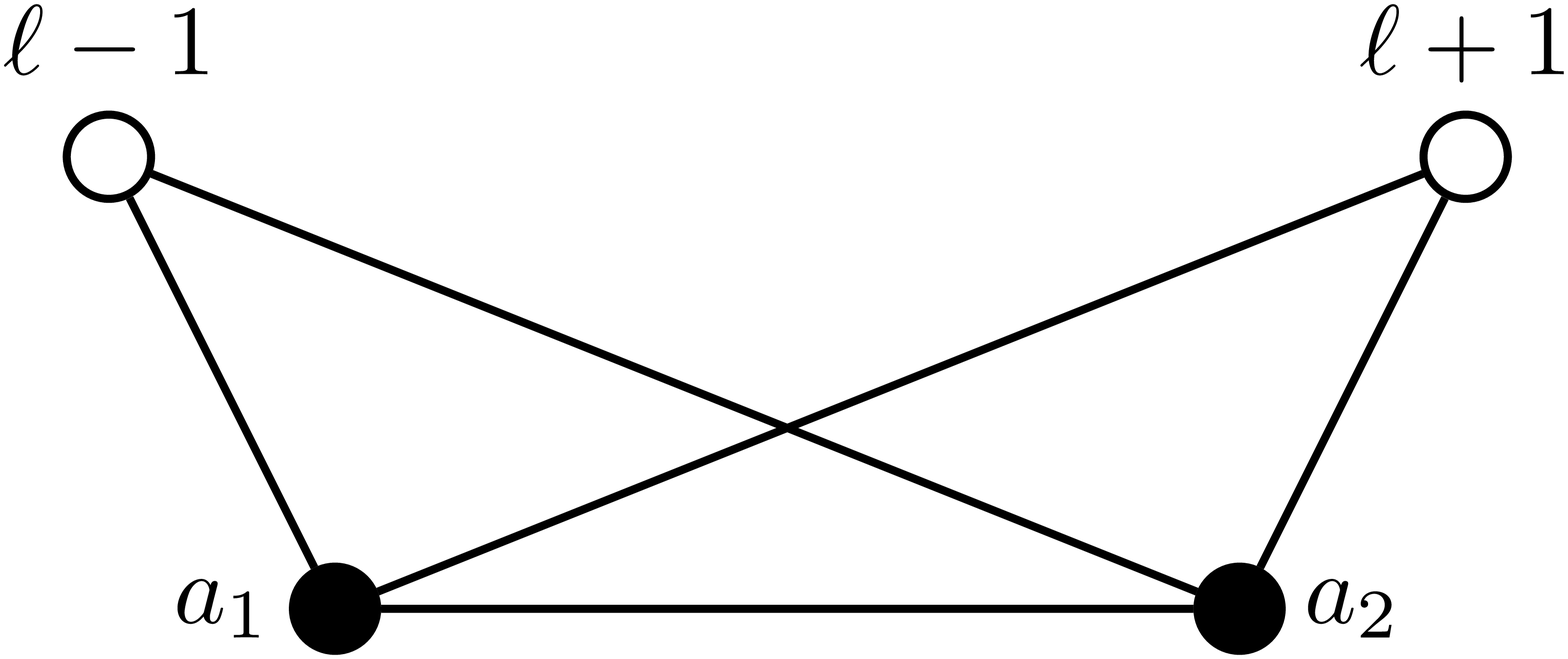} & \includegraphics[width=0.3\textwidth]{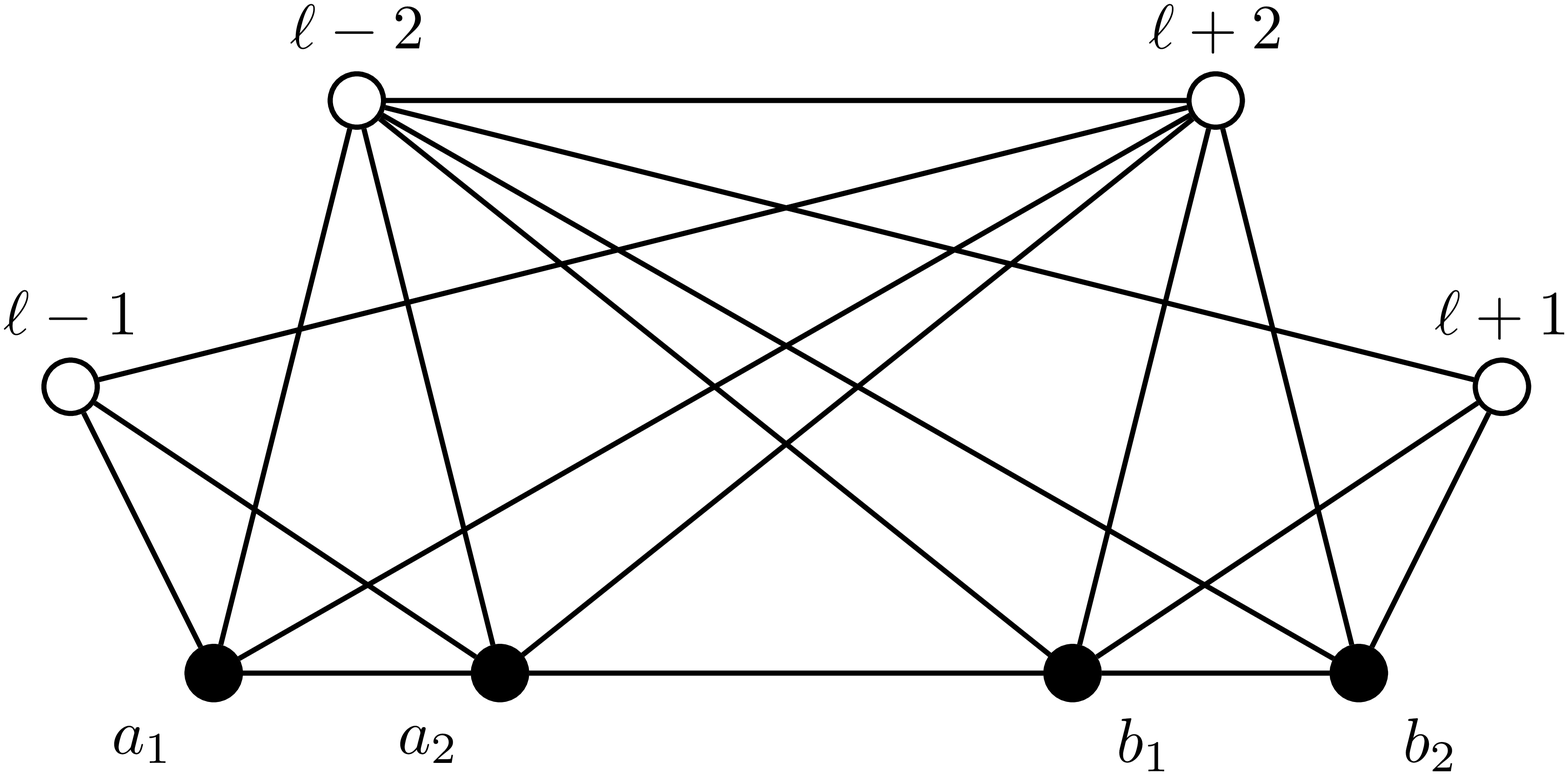} & \includegraphics[width=0.3\textwidth]{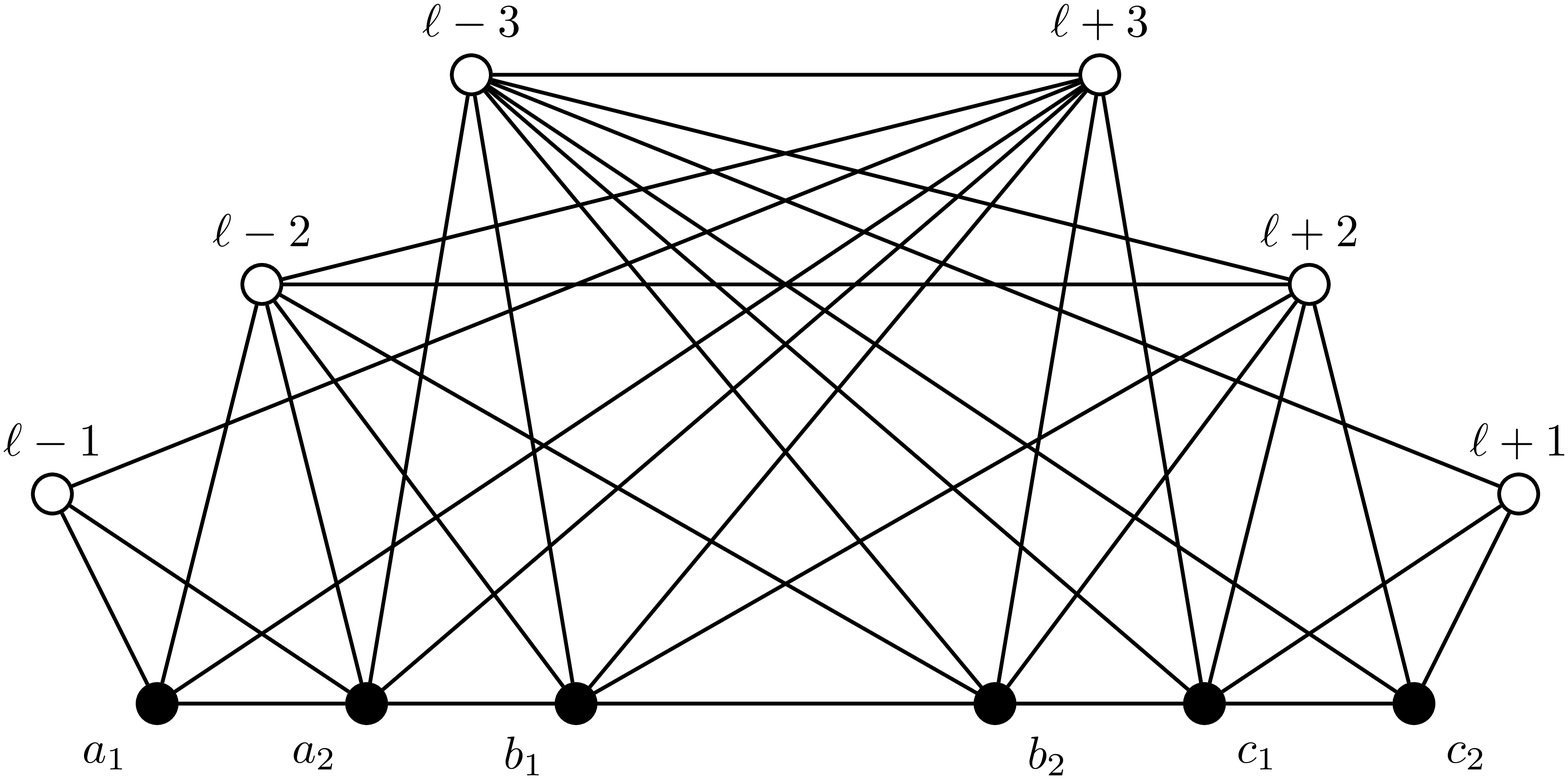}
\end{array}$
\end{center}

Using the formula given in Proposition \ref{generalsumformula}, we can compute $\hstar{\Delta_{9,2}^{stab(2)}}$.  
By Proposition \ref{sumformula} we know that
	$$\hstar{\Delta_{9,2}^{stab(3)}}=1+9x+27x^2+30x^3+9x^4.$$
Hence, it remains to compute the independence polynomials for $G_{9,2,\ell}$ for each $\ell\in[9]$.  
These polynomials are identified in the following table.
\begin{center}
	\includegraphics[width=0.5\textwidth]{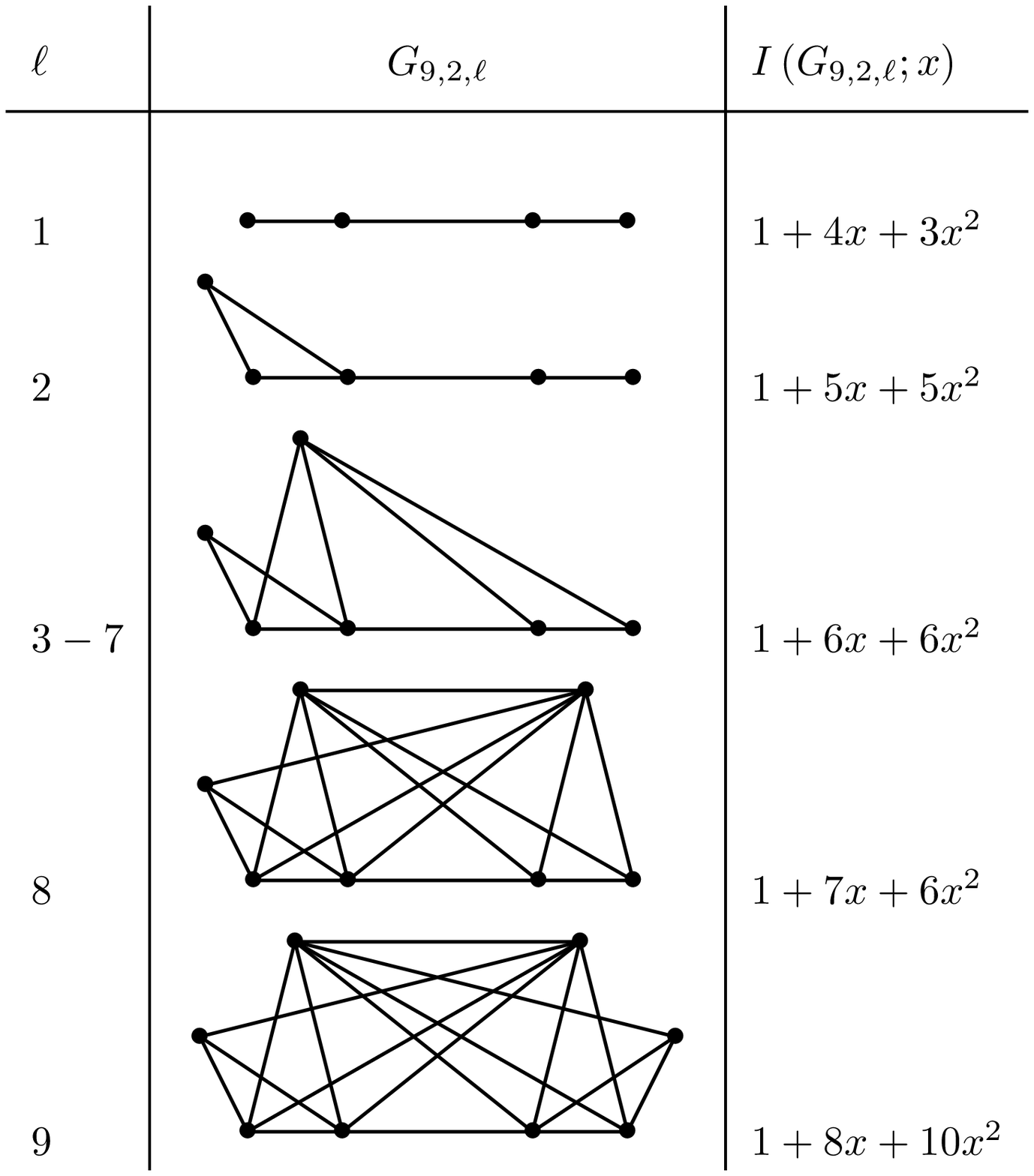}
\end{center}
\end{example}
It follows that 
\begin{equation*}
\begin{split}
	\hstar{\Delta_{9,2}^{stab(2)}}&=\hstar{\Delta_{9,2}^{stab(3)}}+x(9+54x+54x^2),\\
		&=1+18x+81x^2+84x^3+9x^4.\\
\end{split}
\end{equation*}
Similarly, we can also compute that

\begin{equation*}
\begin{split}
\hstar{\Delta_{7,2}^{stab(1)}}&=1+14x+35x^2+7x^3,\\
\hstar{\Delta_{9,2}^{stab(2)}}&=1+18x+81x^2+84x^3+9x^4, \mbox{ and}\\
\hstar{\Delta_{11,2}^{stab(3)}}&=1+22x+143x^2+297x^3+165x^4+11x^5.\\
\end{split}
\end{equation*}

%%%%%%%%%%%%%%%%%%%%%%%%%%%%%%%%%%%%%%%%%%%%%%%%%%%%%%%%%%
%%%%%%%%%%%%%%%%%%%%%%%%%%%%%%%%%%%%%%%%%%%%%%%%%%%%%%%%%%

\subsection{The $h^*$-polynomial of $\Delta_{n,2}^{stab\left(\floorntwo-1\right)}$ and CR Mappings of Lens Spaces.}\label{lens}
Fix odd $n>2$ and $r:=\floorntwo-1$.  
We begin this subsection by demonstrating that the $h^*$-polynomial of $\rstablentwo$ is the independence polynomial of the cycle on $n$ vertices.
As a corollary, we show that the $h^*$-polynomial of the relative interior of $\rstablentwo$ is a univariate specialization of an important class of polynomials in CR geometry.
%Finally, for a general $r$, we show that the $h^*$-polynomial of the relative interior of $\rstablentwo$ is a univariate specialization of the squared Euclidean norm function of polynomial maps that induce a smooth immersions of the Lens Space $L(n,2)$ into $\C^{n+1}$, thereby providing a partial generalization of the $r=\floorntwo-1$ case.

\begin{thm}\label{lucasthm}
Fix odd $n>2$ and $r=\floorntwo-1$.  Then 
	$$\hstar{\rstablentwo}=I(C_n;x)=L_n(x).$$
\end{thm}

To prove Theorem \ref{lucasthm} we first need a lemma that relates Lucas polynomials and Fibonacci polynomials.

\begin{lem}\label{lucasandfib}
For $n\geq 2$, the $n^{th}$ Lucas polynomial can be computed as
	$$L_n(x)=F_n(x)+xF_{n-2}(x).$$
\end{lem}

\begin{proof}
In \cite{levit2} it is noted that the diagonals of Pascal's Triangle given by 
	$${n\choose 0},{n-1\choose 1},{n-2\choose 2},\ldots,{\ceilingntwo\choose \floorntwo}$$
are precisely the coefficients of the $n^{th}$ Fibonnaci polynomial
	$$F_n(x)={n\choose 0}+{n-1\choose 1}x+{n-2\choose 2}x^2+\cdots+{\ceilingntwo\choose \floorntwo}x^{\floorntwo}.$$
In Appendix A of \cite{kappraff} it is noted that the $n^{th}$ Lucas polynomial is given by the same diagonal in the modified Pascal's Triangle

\begin{center}
\begin{tabular}{ccccccccccc}
		&		&		&		&		&	2\\\noalign{\smallskip\smallskip}
		&		&		&		&	1	&		&	2\\\noalign{\smallskip\smallskip}
		&		&		&	1	&		&	3	&		&	2\\\noalign{\smallskip\smallskip}
		&		&	1	&		&	4	&		&	5	&		&	2\\\noalign{\smallskip\smallskip}
		&	1	&		&	5	&		&	9	&		&	7	&		&	2\\\noalign{\smallskip\smallskip}
	1	&		&	6	&		&	14	&		&	16	&		&	9	&		&	2\\\noalign{\smallskip\smallskip}
\end{tabular}
\end{center}
For convenience, we refer to this triangle as Lucas' Triangle.  
One way to produce Lucas' Triangle is to write the $2's$ on the right boundary as $2=1_b+1_g$.  
In this way, we have a blue $1$ and and green $1$ summing to give $2$.  
Imagine that the $1's$ on the left boundary are also blue $1's$.  
Now as we fill in the interior of the triangle using the standard Pascalian recursion write each entry as the sum of the blue $1's$ plus the sum of the green $1's$.  This yields
\begin{center}
\begin{tabular}{cccccccccc}
		&		&		&		&		&	$1_b+1_g$\\\noalign{\smallskip}
		&		&		&		&	$1_b$	&		&	$1_b+1_g$\\\noalign{\smallskip}
		&		&		&	$1_b$	&		&	$(2)_b+(1)_g$	&		&	$1_b+1_g$\\\noalign{\smallskip}
		&		&	$1_b	$&		&	$(3)_b+(1)_g$	&		&	$(3)_b+(2)_g$	&		&	$1_b+1_g$\\\noalign{\smallskip}
	$1_b$	&		&		&	$(4)_b+(1)_g$	&		&	$(6)_b+(3)_g$	&		&	$(4)_b+(3)_g$	&		&	$1_b+1_g$\\\noalign{\smallskip}
%	$1_b$	&		&	$(5)+(1)$	&		&	$(10)+(4)$	&		&	$(10)+(6)$	&		&	$(5)+(4)$	&		&	$1_b+1_g$\\\noalign{\smallskip}
\end{tabular}
\end{center}

With this decomposition we see that Lucas' Triangle can be produced by a term-by-term sum of Pascal's Triangle and the Pascal-like triangle
\begin{center}
\begin{tabular}{ccccccccccc}
		&		&		&		&		&	1\\\noalign{\smallskip\smallskip}
		&		&		&		&	0	&		&	1\\\noalign{\smallskip\smallskip}
		&		&		&	0	&		&	1	&		&	1\\\noalign{\smallskip\smallskip}
		&		&	0	&		&	1	&		&	2	&		&	1\\\noalign{\smallskip\smallskip}
		&	0	&		&	1	&		&	3	&		&	3	&		&	1\\\noalign{\smallskip\smallskip}
	0	&		&	1	&		&	4	&		&	6	&		&	4	&		&	1\\\noalign{\smallskip\smallskip}
\end{tabular}
\end{center}
From this it is then easy to deduce  the identity
	$$L_n(x)=F_n(x)+xF_{n-2}.$$
\end{proof}

With this lemma in hand, we are ready to prove Theorem \ref{lucasthm}.

\subsubsection{Proof of Theorem \ref{lucasthm}.}
\label{proofoflucasthm}
Recall that for a fixed $\ell\in[n]$ the number of simplices in \\
$\max\trintwor\backslash\max\trintworplus$ with maximum $r$-adjacent vertex $\adj{\ell}$ with unique minimal new face of dimension $i$ is the number of ways to construct a lattice path in the $\ell$-region that uses precisely $i$ accessible lattice points on the lines $y=x$ and $y=x-3$.
Recall however, that in the $\ell$-region the lattice points $(0,0)$ and $(n-r,r)$ correspond to the same vertex.  
In particular, we can think of the $\ell$-region as repeating itself in the region translated right $n-r$ and up $r$.  
Reflecting this translated copy of the $\ell$-region about the line $y=x-3$ results in a strip of height $n$ in the region between $y=x$ and $y=x-3$ containing a lattice path, corresponding to $\lambda^\star$, that only touches the boundary diagonals at lattice points labeled by $\ell$.  
Flip the corner of this path at the lattice point labeled by $\ell$ (also making the corresponding flips at the top and bottom of the diagram), and label the two points on the opposite boundary line that are inhibited by $\ell$ as $\ell-(r+1)$ and $\ell+(r+1)$.  
This results in a diagram with $(t,t)$ labeled by $\ell+t$ and $(t,t-3)$ labeled by $\ell-r+t$ for $t\in[n]$.
The diagram below represents these manipulations of the $\ell$-region for $n=9$.

\begin{center}
	\includegraphics[width=1.0\textwidth]{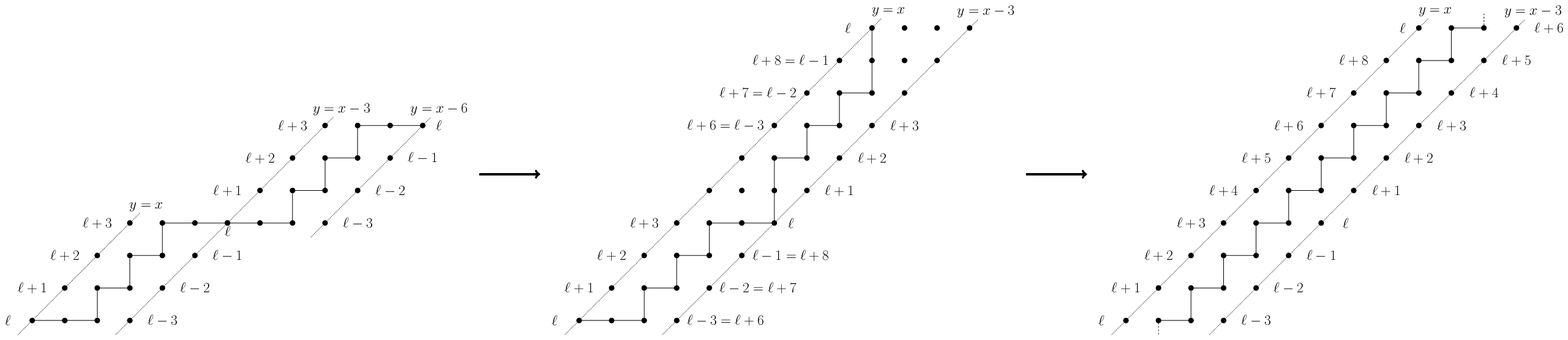}
\end{center}

For this region, the inhibition diagram for the spots 
	$$X=\{1,2,3,\ldots,n\}$$
is exactly a cycle on $n$ vertices labeled as 

\begin{center}
	\includegraphics[width=0.4\textwidth]{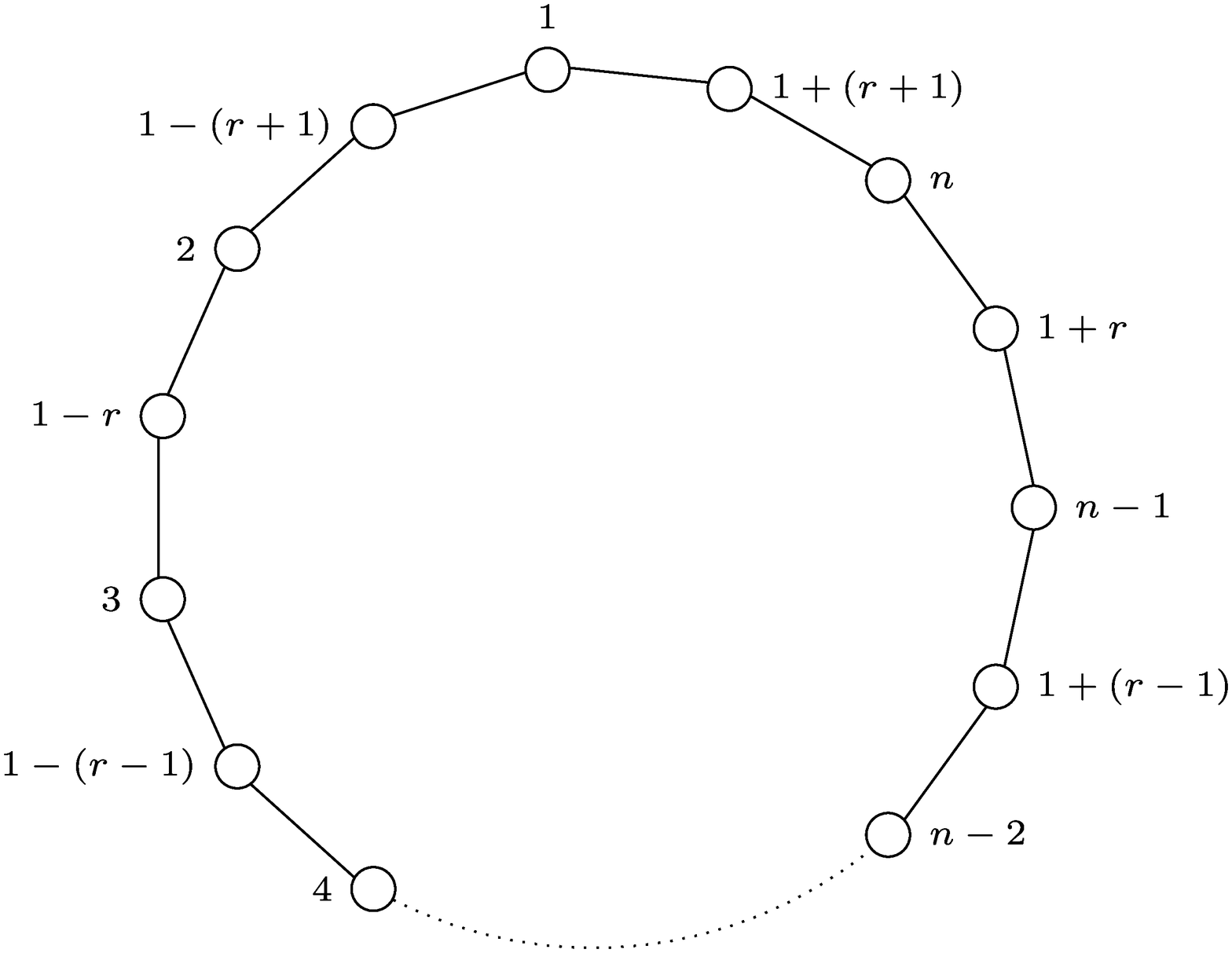}
\end{center}

\noindent We remark that this labeling of $C_n$ is precisely the underlying graph of the edge polytope $\Delta_{n,2}^{stab\left(\floorntwo\right)}$, the next smallest $r$-stable hypersimplex in the chain.

Notice now that for a simplex with maximum $r$-adjacent vertex $\adj{\ell}$ and unique minimal new face of dimension $i$ we have a unique independent set in this labeled copy of $C_n$ of cardinality $i+1$ with maximal vertex label being $\ell$.  
Conversely, if we pick an arbitrary independent $(i+1)$-set in this labeled copy of $C_n$ it has a unique maximally labeled element, say the vertex labeled by $\ell$.  
In the region for which this inhibition diagram arises we may then flip the corners of the lattice path corresponding to $\lambda^\star$ such that it touches the boundary diagonals $y=x$ and $y=x-3$ precisely at the lattice points labeled by the elements of the independent set.  
Then translate the diagram so that $\ell$ is labeling the origin (and the first move is an East move).  
This lattice path in this $\ell$-region gives the corresponding simplex in $\max\trintwor\backslash\max\trintworplus$.  
This establishes a bijection between the independent $(i+1)$-sets in this labeled copy of $C_n$ and the simplices in $\max\trintwor\backslash\max\trintworplus$ with unique minimal new face of dimension $i$.
Hence,
	$$\hstar{\rstablentwo}=I(C_n;x).$$

Finally, to see that this polynomial is also the $n^{th}$ Lucas polynomial we utilize a result of Arocha \cite{arocha} which states that
	$$I(C_n;x)=F_{n-1}(x)+2xF_{n-2}(x).$$
From this fact, and the identity from Lemma \ref{lucasandfib}, we see that
	\begin{equation*}
	\begin{split}
	\hstar{\rstablentwo}&=I(C_n;x),\\
		&=F_{n-1}(x)+2xF_{n-2}(x),\\
		&=(F_{n-1}(x)+xF_{n-2}(x))+xF_{n-2}(x),\\
		&=F_n(x)+xF_{n-2}(x),\\
		&=L_n(x).\\
	\end{split}
	\end{equation*}
This completes the proof of Theorem \ref{lucasthm}.
\qed

\bigskip
\bigskip

We now describe a connection between our polytopes and CR geometry.  
CR geometry is a fascinating field of study that examines properties of real hypersurfaces as submanifolds of $\C^n$ and their intrinsic complex structure induced by the ambient space $\C^n$.  
The interested reader should refer to \cite{dangelo1} for a nice introduction to this theory.  
In \cite[Chapter 5]{dangelo1}, D'Angelo describes the theory of proper holomorphic mappings between balls that are invariant under subgroups of the unitary group $U(n)$.  
These maps are particularly interesting as they induce maps from spherical space forms into spheres.  
In \cite{dangelo3}, D'Angelo, Kos, and Riehl define recursively the following collection of polynomials (the former-most author also defines this collection explicity in \cite{dangelo1, dangelo2}).
Let 
	$$g_0(x,y)=x, \qquad g_1(x,y)=x^3+3xy, $$
and 
	$$g_n(x,y)=(x^2+2y)g_{n-1}(x,y)-y^2g_{n-2}(x,y),$$
for $n\geq2$.
Then set 
	$$p_n(x,y)=g_n(x,y)+y^{2n+1}.$$
For odd $n>2$ consider the group $\Gamma(n,2)$ of $2\times 2$ complex matrices of the form
\[ \left( \begin{array}{cc}
\gamma & 0 \\
0 & \gamma^2 \\
\end{array} \right)^k \]
where $\gamma$ is a primitive $n^{th}$ root of unity and $k=0,1,2,\ldots,n-1$.  
Recall that the Lens space, $L(n,2)$ is defined as the quotient $L(n,2)=S^3/\Gamma(n,2)$.  
In \cite{dangelo1,dangelo2} it is shown that the polynomials $p_n(x,y)$ correspond to proper holomorphic monomial maps between spheres (given by their monomial components) that are invariant under $\Gamma(n,2)$.  
In \cite{dangelo3} it is shown that these polynomials are of highest degree with respect to this property.
These polynomials share the following relationship with $\rstablentwo$ via Ehrhart-MacDonald reciprocity.

\begin{thm}\label{CRgeometrythm}
Fix odd $n>2$ and let $r=\floorntwo-1$.  Then 
	$$p_{r+1}(x,x)=\hstar{\rstablentwoint}+x^n,$$
where $\rstablentwoint$ denotes the relative interior of $\rstablentwo$.
Hence, the $h^\ast$-polynomial of the relative interior of $\rstablentwo$ (plus an $x^n$ term) is a univariate evaluation of the squared Euclidean norm function of a monomial CR mapping of the Lens space, $L(n,2)$, into the unit sphere of complex dimension $r+3$.
\end{thm}

\begin{proof}
Consider the Lucas sequence defined by the recurrence
	$$\el_0(x,y)=2, \quad \el_1(x,y)=x, \quad \mbox{ and } \quad \el_n(x,y)=x\el_{n-1}(x,y)+y\el_{n-2}(x,y)$$
for $n\geq 2$.  
Comparing this sequence to the Lucas polynomials $\{L_n\}$ we may deduce that
	$$\el_n(x,y)=x^nL_n\left(\frac{y}{x^2}\right).$$
By applying the recurrence, we see that
	$$\el_{n+2}(x,y)=\el_n(x,y)(x^2+2y)-y^2\el_{n-2}(x,y).$$
It then follows that the odd terms of the sequence $\{\el_n\}$ are precisely the terms of the sequence $\{g_n(x,y)\}$.  
That is, if $n=2m+1$ then 
	$$\el_n(x,y)=g_m(x,y).$$
Recall now that $n=2r+3=2(r+1)+1$.
Applying these facts, together with Ehrhart-MacDonald reciprocity, we see that
	\begin{equation*}
	\begin{split}
	p_{r+1}(x,x)&=g_{r+1}(x,x)+x^n,\\
		&=\el_n(x,x)+x^n,\\
		&=x^nL_n\left(\frac{x}{x^2}\right)+x^n,\\
		&=x^n\hstarint{\rstablentwo}+x^n,\\
		&=\hstar{\rstablentwoint}+x^n.\\
	\end{split}
	\end{equation*}
The fact that $p_{r+1}(x,y)$ is the squared Euclidean norm function of a monomial CR mapping of the Lens space, $L(n,2)$, into a sphere of complex dimension $r+3$ is proven via the discussion in \cite[pp.171-174]{dangelo1}.
\end{proof}

More generally, for $r\leq \floorntwo-1$, the $h^\ast$-polynomial of the relative interior of $\rstablentwo$ is a univariate specialization of the squared Euclidean norm function of certain polynomial maps that induce smooth immersions of the Lens space $L(n,2)$ into $\C^{\floorntwo+1}$. 
However, it is unclear that these maps are mapping $L(n,2)$ into a hypersurface with interesting structure.  
For this reason we pose the following question.

\begin{quest}\label{lensspacequestion}
For $r<\floorntwo-1$, does the $h^\ast$-polynomial of the relative interior of $\rstablentwo$ arise as a univariate specialization of a polynomial corresponding to maps between interesting manifolds?
\end{quest}

%%%%%%%%%%%%%%%%%%%%%%%%%%%%%%%%%%%%%%%%%%%%%%%%%%%%%%%%%%%%%%%%%%%%%%%%%%%%%%%%%%%
%%%%%%%%%%%%%%%%%%%%%%%%%%%%%%%%%%%%%%%%%%%%%%%%%%%%%%%%%%%%%%%%%%%%%%%%%%%%%%%%%%%
%%%UNIMODALITY STUFF %%%%%%%%%%%%%%%%%%%%%%%%%%%%%%%%%%%%%%%%%%%%%%%%%%%%%%%%%%%%%%%%%%%%
%%%%%%%%%%%%%%%%%%%%%%%%%%%%%%%%%%%%%%%%%%%%%%%%%%%%%%%%%%%%%%%%%%%%%%%%%%%%%%%%%%%
%%%%%%%%%%%%%%%%%%%%%%%%%%%%%%%%%%%%%%%%%%%%%%%%%%%%%%%%%%%%%%%%%%%%%%%%%%%%%%%%%%%

\subsection{Unimodality}
A consequence of work by Katzman \cite{kat} is that the $h^\ast$-polynomial of $\hyperntwo$ is unimodal.  
It appears that this is also true for the $r$-stable hypersimplices within.  
In this subsection, we utilize the shelling and our computations of $h^\ast$-polynomials to show that this observed unimodality does in fact hold in some specific cases.
We begin with two quick corollaries, one to Theorem \ref{lucasthm} and the other to Theorem~\ref{even gorenstein hstarpolynomials theorem}.

\begin{cor}
\label{logconcave}
Fix odd $n>2$ and let $r=\floorntwo-1$.  The $h^\ast$-polynomial of $\rstablentwo$ is log-concave and hence unimodal.
\end{cor}

\begin{proof}
Since the $h^\ast$-polynomial of $\rstablentwo$ is $I(C_n;x)$ and $C_n$ is a claw-free graph then by \cite{hami} it is log-concave and consequently unimodal.
\end{proof}

\begin{cor}
\label{even gorenstein unimodality corollary}
Let $n=2r+2$.  The $h^\ast$-polynomial of $\rstablentwo$ is log-concave and unimodal.  
\end{cor}

\begin{proof}
This corollary is immediate for the result of Theorem~\ref{even gorenstein hstarpolynomials theorem}, which shows this $h^\ast$-polynomial is the generating polynomial for the binomial coefficiencts ${r+1\choose m}$.  
\end{proof}

We also note that the polynomials described in Corollaries~\ref{logconcave} and \ref{even gorenstein unimodality corollary} possess the much stronger property of \emph{real-rootedness}, meaning that all zeros of these $h^\ast$-polynomials are real numbers.  
The following results utilize the formula for the $h^\ast$-polynomial of $\Delta_{n,2}$ given by Katzman in \cite{kat}.  
Using this formula, we subtractively compute formulas for the $h^\ast$-polynomials of $\twostablentwo$ and $\threestablentwo$ by ``undoing" our shelling.

\begin{cor}\label{twostableunimodal}
The $h^\ast$-vector of $\twostablentwo$ is unimodal for $n$ odd.
\end{cor}

\begin{proof}
Consider the completion of the shelling of $\nabla_{n,2}^2$ to a shelling of $\nabla_{n,2}$.  Here, $r=1$, and so the simplices $\omega_{\ell,\lambda}\in\max\nabla_{n,2}\backslash\max\nabla_{n,2}^2$ are labeled with compositions of length $n-2$.  The composition $\lambda$ is a composition of $1$ that must satisfy equation (\ref{lambdabounds}).  We now determine which compositions are admissible for a fixed $\adjvar{1}{\ell}=\max v_{(\omega)}\cap\Adjvar{1}{n}$.  

For $\ell\in\{2,3,4,\ldots,n-1\}$ the composition $(1,0,0,\ldots,0)$ does not label a simplex with $\adjvar{1}{\ell}=\max v_{(\omega)}\cap\Adjvar{1}{n}$, since such a composition would necessarily have $\max v_{(\omega)}\cap\Adjvar{1}{n}=\adjvar{1}{\ell+1}$.  This is depicted in Figure \ref{fig:unimodaltwostable1}.  On the other hand, each other composition of $1$ into $n-2$ parts does correspond to a simplex with $\adjvar{1}{\ell}= \max v_{(\omega)}\cap\Adjvar{1}{n}$.  By considering the associated lattice paths in Figure \ref{fig:unimodaltwostable2} it is easy to see that the simplex $\omega_{\ell,\lambda^\star}$ has $\#(G_\omega)=1$, and the other $n-4$ simplices $\omega_{\ell,\lambda}$ with $\lambda=(0,0,\ldots,0,1,0,\ldots,0)$ have $\#(G_\omega)=2$.  

\begin{figure}
	\centering
	\includegraphics[width=0.6\textwidth]{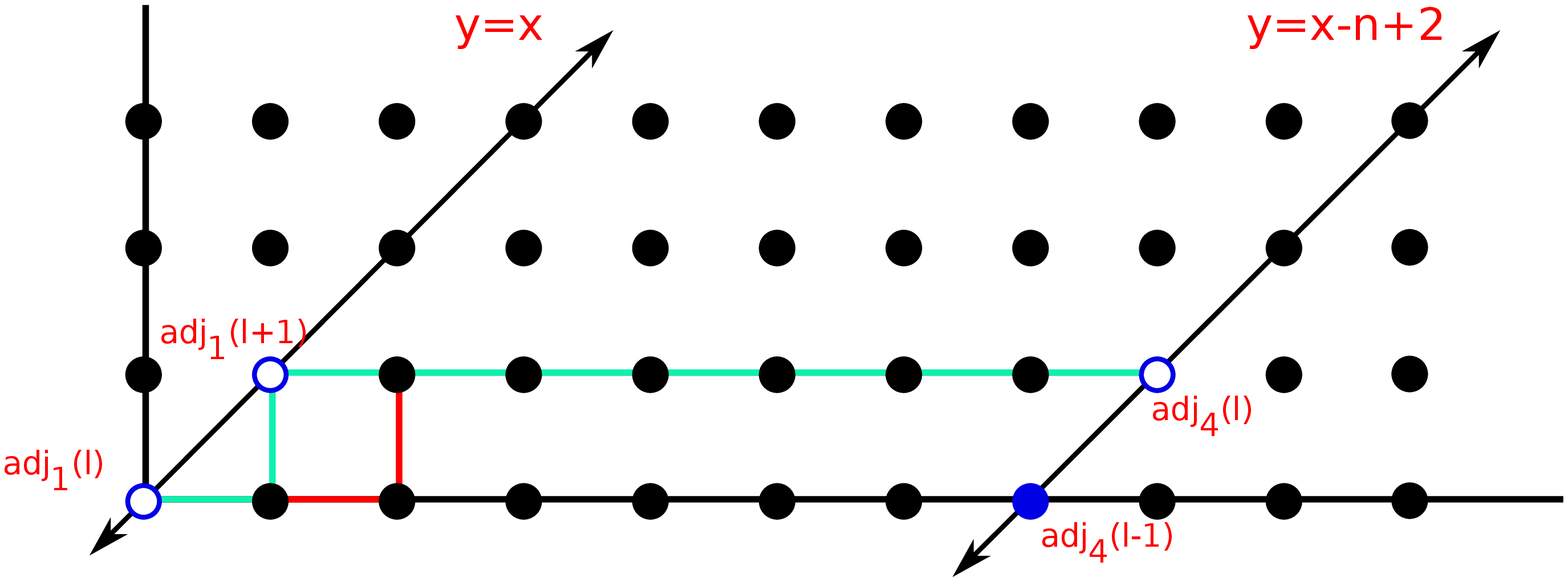}
	\caption{Here we let $n=9$.}
	\label{fig:unimodaltwostable1}
\end{figure}

For $\ell=1$, the compositions $(1,0,0,\ldots,0)$ and $(0,0,\ldots,0,1)$ do not label a simplex with $\adjvar{1}{1}=\max v_{(\omega)}\cap\Adjvar{1}{n}$, since such a simplex necessarily has $\max v_{(\omega)}\cap\Adjvar{1}{n}\in\{2,n\}$.  This is depicted in Figure \ref{fig:unimodaltwostable2}.  Again, the simplex $\omega_{1,\lambda^\star}$ has $\#(G_\omega)=1$, and the remaining $n-5$ simplices $\omega_{1,\lambda}$, for $\lambda=(0,0,\ldots,0,1,0,\ldots,0)$ have $\#(G_\omega)=2$.

\begin{figure}
	\centering
	\includegraphics[width=0.6\textwidth]{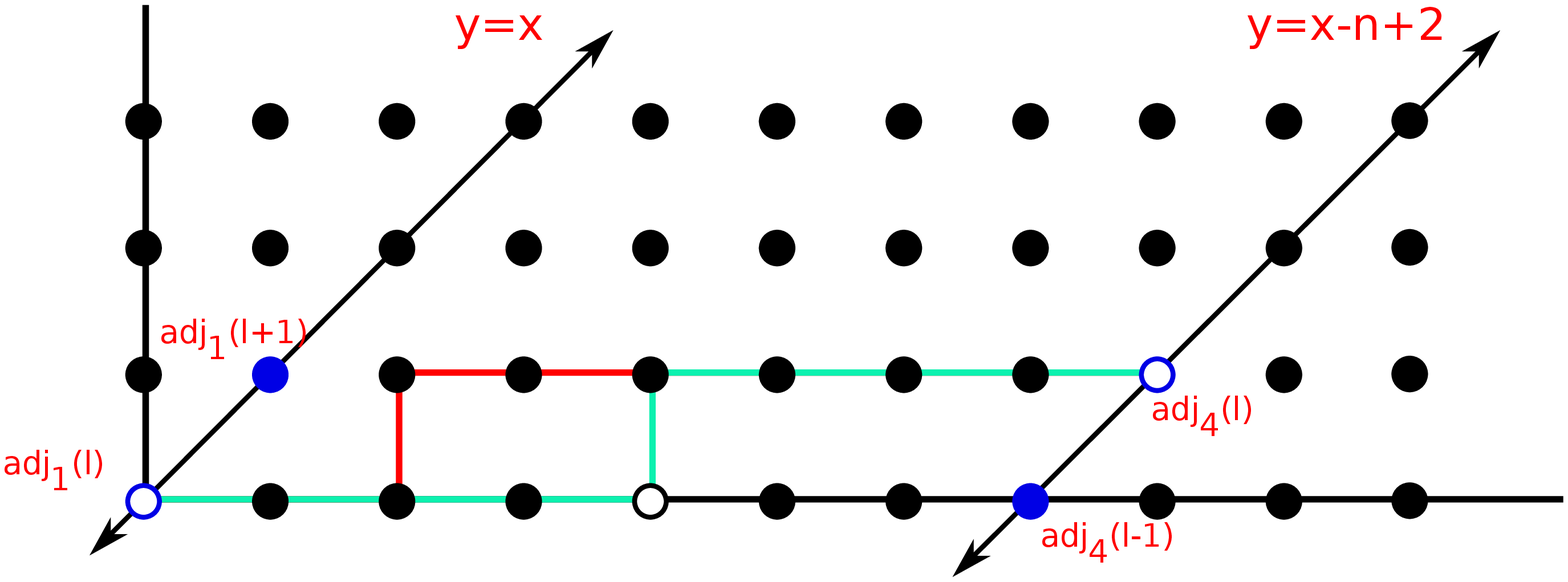}
	\caption{Here we let $n=9$.}
	\label{fig:unimodaltwostable2}
\end{figure}

Finally, for $\ell=n$ both the compositions $(1,0,0,\ldots,0)$ and $(0,0,\ldots,0,1)$ label a simplex with $\adjvar{1}{1}=\max v_{(\omega)}\cap\Adjvar{1}{n}$, since $\adjvar{1}{n}=\max\Adjvar{1}{n}$.  Hence there is one simplex, namely $\omega_{n,\lambda^\star}$, with $\#(G_\omega)=1$, and the remaining $n-3$ simplices have $\#(G_\omega)=2$.  Summarizing this analysis we have $n$ simplices with $\#(G_\omega)=1$, and $n(n-4)$ simplices with $\#(G_\omega)=2$.  

In \cite{kat}, Katzman computed that for $\hyperntwo$
$$h_i^\ast={n \choose 2i}$$
for $i\neq1$, and
$$h_1^\ast={n\choose 2}-n.$$
Thus, since there are $n$ elements in $\Adjvar{1}{n}$ we have that
\begin{equation*}
\begin{split}
(h_1^\ast)^{stab(2)}&=h_1^\ast-n={n\choose 2}-2n,\\
(h_2^\ast)^{stab(2)}&=h_2^\ast-n(n-4)={n\choose 4}-n(n-4), \mbox{ and for $i\neq1,2$}\\
(h_i^\ast)^{stab(2)}&=h_i^\ast.\\
\end{split}
\end{equation*}
It is then easy to verify that $(h_1^\ast)^{stab(2)}-(h_0^\ast)^{stab(2)}\geq0$ and $(h_2^\ast)^{stab(2)}-(h_1^\ast)^{stab(2)}\geq0$ for all $n\geq0$.  As well, $(h_3^\ast)^{stab(2)}-(h_2^\ast)^{stab(2)}\geq0$ for all $n\neq6,7,8$.  But this is fine since the $h^\ast$-vector for $\Delta_{7,2}^{stab(2)}$ is
	$$h^\ast\left(\Delta_{7,2}^{stab(2)}\right)=(1,7,14,7,0,0).$$
We also remark that
	$$h^\ast\left(\Delta_{6,2}^{stab(2)}\right)=(1,3,3,1,0), \mbox{ and}$$
	$$h^\ast\left(\Delta_{8,2}^{stab(2)}\right)=(1,12,38,28,1,0,0).$$
\end{proof}

Notice that the result given by this subtractive formula for $\hstar{\Delta_{9,2}^{stab(2)}}$ agrees with the result via independence polynomials computed in Example \ref{inhibitiondiagrams}.
It is possible to apply the same strategy used in the proof of Corollary \ref{twostableunimodal} to show that the $h^\ast$-vector of $\threestablentwo$ is unimodal.  
In short, we count the lattice paths corresponding to simplices in the set $\max\nabla_{n,2}^2\backslash\max\nabla_{n,2}^3$ with unique minimal new face of cardinality $i=1,2,3$ for each choice of maximal $\adj{\ell}$, $\ell\in[n]$.  
We then subtract these values from the corresponding coefficients in the $h^\ast$-vector of $\twostablentwo$, and check that the unimodality condition is satisfied for the resulting $h^\ast$-vector.  
However, the details of this computation are quite unpleasant, so we omit them.

\begin{cor}\label{threestableunimodal}
Let $n$ be odd.  The $h^\ast$-vector of $\threestablentwo$ is given by

\begin{equation*}
\begin{split}
(h_1^\ast)^{stab(3)}&={n\choose 2}-3n,\\
(h_2^\ast)^{stab(3)}&={n\choose 4}-\frac{1}{2}(n(7n-55)+94),\\
(h_3^\ast)^{stab(3)}&={n\choose6}-\frac{1}{2}(n^3-13n^2+40n+16), \mbox{ and for $i\neq1,2,3$}\\
(h_i^\ast)^{stab(3)}&={n \choose 2i}.\\
\end{split}
\end{equation*}
\noindent Moreover, $h^\ast\left(\threestablentwo\right)$ is unimodal.
\end{cor}

We end this section with two conjectures that arise from these observations. 
\begin{conj}
\label{second hypersimplices conjecture}
The $h^\ast$-polynomials of the $r$-stable second hypersimplices $\rstablentwo$ are unimodal.  
\end{conj}
A natural first-step to validating Conjecture~\ref{second hypersimplices conjecture} would be to prove the following.
\begin{conj}
\label{independence polynomials conjecture}
The independence polynomials $I(G_{n,r,\ell};x)$ are unimodal. 
\end{conj}

%%%%%%%%%%%%%%%%%%%%%%%%%%%%%%%%%%%%%%%%%%%%%%%%%%%%%%%%%%%%%%%%%%%%%%%%%%%%%%%%%%%
%%%%%%%%%%%%%%%%%%%%%%%%%%%%%%%%%%%%%%%%%%%%%%%%%%%%%%%%%%%%%%%%%%%%%%%%%%%%%%%%%%%
\section{The $r$-stable $(n,k)$-hypersimplices for $k>2$.}
\label{extending results}

In \cite{lam}, Lam and Postnikov establish that a minimal circuit $(\omega)$ in $G_{n,k}$ corresponds to a collection of $n$ permutations in $\sym_n$ that form an equivalence class modulo cyclic shifts
$$
\omega_1\omega_2\cdots\omega_n\sim\omega_n\omega_1\omega_2\cdots\omega_{n-1},
$$
and each permutation $\omega$ in this equivalence class corresponds to a different choice of initial vertex of the circuit $(\omega)$, which we denote $\epsilon(\omega)$.  
Here, the permutation $\omega^{-1}$ encodes the location of the $1$'s in the $(0,1)$-vector $\epsilon(\omega)$ via its descents, and the permutation $\omega$ provides an algorithm to reproduce all vertices of the minimal circuit given this initial vertex.  
Often times, it is easier to work with the simplices of $\trink$ when they are encoded as permutations.  
In this section, we will rephrase the shelling order $<$ on $\max\trintwo$ in terms of permutations, and present an order on $\max\trink$ for $k>2$ which the authors conjecture generalizes the shelling order.  
We first develop some terminology.

For a permutation $\omega=\omega_1\omega_2\cdots\omega_n\in\sym_n$ let
$$
\Des(\omega):=
\{i\in[n] :
\mbox{ $\omega_i>\omega_{i+1}$ or $i=n$ and $\omega_n>\omega_1$}
\}
$$
be the set of \emph{extended descents of $\sym_n$}.
For a minimal circuit $(\omega)$ in $G_{n,k}$ and a fixed representative $\omega\in\sym_n$ of $(\omega)$ the inverse permutation $\omega^{-1}$  will have exactly $k$ extended descents, and
\begin{equation}
\label{extended descents equation}
\Des(\omega^{-1})=
\{ 
i\in[n] :
\epsilon(\omega)_{i+1}=1
\}.
\end{equation}
Let $\omega^i$ denote the inverse permutation of the representative of $(\omega)$ given by letting the move $\epsilon\overset{n}\rightarrow \epsilon^\prime$ in $(\omega)$ be the $i^{th}$ move in the circuit.   
Then the minimal circuit $(\omega)$ is encoded via the $n$-tuple of permutations
$$
\W(\omega) :=
(\omega^1,\omega^2,\ldots,\omega^n)\in(\sym_n)^{\times n}.
$$
The $n^{th}$ coordinate of $\W(\omega)$ is the typical encoding of $(\omega)$ used in \cite{stan2} and \cite{lam}. 
The encoding of $(\omega)$ as $\omega^n$ can be thought of as encoding the whole circuit via an initial condition (the initial vertex) and an algorithm with which to recover the circuit given the initial condition (the permutation).  
The $n$-tuple $\W(\omega)$ simply encodes all possible pairs of initial conditions and their algorithms.  

Generalizing our previous definition, an \emph{$r$-adjacent vertex} $\epsilon=(\epsilon_1,\epsilon_2,\ldots,\epsilon_n)$ of $\hypernk$ is any vertex of $\hypernk$  for which there exists an $i\in[n]$ such that $\epsilon_i=\epsilon_{i+r}$ where we consider the indices modulo $n$.  
Similarly, we define an \emph{$r$-adjacent extended descent} to be an ordered pair $(i,i+r)\in\Des(\omega)\times\Des(\omega)$, 
and we denote the set of all such pairs by $\Des_r(\omega)$.  
We then consider the map defined coordinate-wise:
\begin{equation*}
\begin{split}
\Phi&: (\sym_n)^{\times n}\longrightarrow \ZZ^n_{\geq0};\\
\Phi&: \W_i\longmapsto \#\Des_r(\W_i),\\
\end{split}
\end{equation*}
applied to the collection
$$
\sym_{n,k}^r:=
\{
\W\in(\sym_n)^{\times n} :
\W=\W(\omega), \omega\in\max\trinkr\backslash\max\trink^{r+1}
\}.
$$

Let $k=2$ and recall that we ordered the simplices $\omega_{\ell,\lambda}\in\max\trintwor\backslash\max\trintwo^{r+1}$ via the colexicographic ordering on the sets $W_{\ell,s}$ refined by the colexicographic ordering on the compositions $\lambda$.  
With this new language, $\sum_{i=1}^n\Phi(\W(\omega_{\ell,\lambda}))_i=s$, the coordinate $\Phi(\W(\omega_{\ell,\lambda}))_\ell$ is the right-most nonzero entry of $\Phi(\W(\omega_{\ell,\lambda}))$, and the coordinate $\W(\omega_{\ell,\lambda})_\ell$ of $\W(\omega_{\ell,\lambda})$ has initial vertex $\adj{\ell}$.  
Thus, the colexicographic ordering on the sets $W_{\ell,s}$ corresponds to the graded colexicographic ordering on the points in the set $\Phi(\sym_{n,2}^r)\subset\ZZ^n_{\geq0}$.  
To refine this partial order on $\sym_{n,2}^r$ in the same fashion as the colexicographic order on the compositions $\lambda$, we examine the permutation structure of the coordinate $\W(\omega_{\ell,\lambda})_\ell$ of $\W(\omega_{\ell,\lambda})$.  
\begin{lem}
\label{permutation lemma}
Let $\omega_{\ell,\lambda},\omega_{\ell,\lambda^\prime}\in W_{\ell,s}$ with $\lambda<_{\colex}\lambda^\prime$.  
Then
$$
\Des_r(\W(\omega_{\ell,\lambda})_\ell)
=\Des_r(\W(\omega_{\ell,\lambda^\prime})_\ell)
=\left\{(\ell-1,\ell+r-1)\right\},
$$
and 
$$
(\W(\omega_{\ell,\lambda})_\ell)_\ell(\W(\omega_{\ell,\lambda})_\ell)_{\ell+1}\cdots(\W(\omega_{\ell,\lambda})_\ell)_{\ell+r-1}
<_{\colex}
(\W(\omega_{\ell,\lambda^\prime})_\ell)_\ell(\W(\omega_{\ell,\lambda^\prime})_\ell)_{\ell+1}\cdots(\W(\omega_{\ell,\lambda^\prime})_\ell)_{\ell+r-1} 
$$
as subwords of $\W(\omega_{\ell,\lambda})_\ell$ and $\W(\omega_{\ell,\lambda^\prime})_\ell$ with indices taken modulo $n$.
\end{lem}

\begin{proof}
Since $\omega_{\ell,\lambda}$ and $\omega_{\ell,\lambda^\prime}$ have the same $r$-adjacent vertex, namely $\adj{\ell}$, then 
$$
\Des_r(\W(\omega_{\ell,\lambda})_\ell)
=\Des_r(\W(\omega_{\ell,\lambda^\prime})_\ell)
=\left\{(\ell-1,\ell+r-1)\right\},
$$
This follows from equality~(\ref{extended descents equation}) and the fact that $\epsilon(\W(\omega_{\ell,\lambda})_\ell)=\adj{\ell}$.  
Finally, for $\ell\leq i\leq \ell+r-1$ we have that $(\W(\omega_{\ell,\lambda})_\ell)_i$ is when $1^L$ makes its $i^{th}$ move in the circuit $\omega_{\ell,\lambda}$.  
Since $\lambda<_{\colex}\lambda^\prime$ then if coordinate $i$ is the right-most coordinate in which these two subwords disagree then it must be that
$$
(\W(\omega_{\ell,\lambda})_\ell)_i<(\W(\omega_{\ell,\lambda^\prime})_\ell)_i.
$$
\end{proof}

By Lemma~\ref{permutation lemma}, we see that the order $<$ on $\max\trintwor\backslash\max\trintwo^{r+1}$ is the same as the graded colexicographic order on $\Phi(\sym_{n,2}^r)$ refined by the colexicographic order on the length $r$ subwords of the permutations $\W(\omega_{\ell,\lambda})_\ell$ specified by the $r$-adjacent extended descents.  
In this new language, the natural generalization to an order on $\sym_{n,k}^r$ for $k>2$ goes as follows.  

First, fix a total ordering on the set of all $r$-adjacent vertices $\Adj{n}$, and partially order the elements of $\sym_{n,k}^r$ via the graded colexicographic order on $\Phi(\sym_{n,k}^r)$.  
If two $n$-tuples $\W_1,\W_2\in \sym_{n,k}^r$ satisfy
$$
\Phi(\W_1)=\Phi(\W_2)
$$ 
then they have the same maximal $r$-adjacent vertex, say $\adjr$, and
$$
\epsilon((\W_1)_\ell)=\adjr=\epsilon((\W_2)_\ell)
$$
for some coordinate $\ell\in[n]$.  
The permutations $(\W_1)_\ell$ and $(\W_2)_\ell$  then have the same extended descent set which we order as 
$$
\Des((\W_1)_\ell)=\{i_1<i_2<\cdots<i_m\}
$$
via the standard ordering on $[n]$.  
They also have the same $r$-adjacent extended descents
$$
\{(j_1,j_1+r),(j_2,j_2+r),\ldots,(j_t,j_t+r)\},
$$
We will let 
$$
\{\rho_1<\rho_2<\cdots<\rho_\alpha\}=\Des((\W_1)_\ell)\backslash\{j_1,\ldots,j_t\}, 
$$
and we let $q_i$ denote the element of $\Des((\W_1)_\ell)$ immediately succeeding $\rho_i$ (modulo $n$) in the given order.  
We then refine the partial order on $\sym_{n,k}^r$ via the colexicographic order on the subwords of the permutations $(\W_1)_\ell$ and $(\W_2)_\ell$ determined by the extended descents.  
As in the $k=2$ case, we first consider those subwords determined by the $r$-adjacent extended descents.
The order is $\W_1<\W_2$ if and only if
for all $i=t,t-1,t-2,\ldots, q$ we have that
\begin{equation}
\label{subword equality}
((\W_1)_\ell)_{j_i}\cdots((\W_1)_\ell)_{j_i+r-1}=((\W_2)_\ell)_{j_i}\cdots((\W_2)_\ell)_{j_i+r-1},
\end{equation}
and for $i=q-1$ we have
$$
((\W_1)_\ell)_{j_i}\cdots((\W_1)_\ell)_{j_i+r-1}<((\W_2)_\ell)_{j_i}\cdots((\W_2)_\ell)_{j_i+r-1},
$$
or if the equality~(\ref{subword equality}) holds for all $i\in[t]$ and for all 
$i=\alpha,\alpha-1,\ldots,q$ we have 
\begin{equation}
\label{subword equality 2}
((\W_1)_\ell)_{\rho_i}\cdots((\W_1)_\ell)_{q_i}=((\W_2)_\ell)_{\rho_i}\cdots((\W_2)_\ell)_{q_i},
\end{equation}
and for $i=q-1$ 
$$
((\W_1)_\ell)_{\rho_i}\cdots((\W_1)_\ell)_{q_i}<((\W_2)_\ell)_{\rho_i}\cdots((\W_2)_\ell)_{q_i}.
$$
This provides a total order on the set $\sym_{n,k}^r$, and hence on $\max\trink$, and this order generalizes the shelling order $<$ on $\max\trintwo$.
Thus, we make the following conjecture.
\begin{conj}
\label{shelling order conjecture}
The given order on $\max\trink$ is a shelling order.
\end{conj}

\end{document}